\documentclass[%
 aip,
 jmp,
 amsmath,amssymb,
 reprint,%
]{revtex4-1}

\usepackage{mathrsfs}
\usepackage{amssymb}
 \usepackage{epsfig}
\usepackage{graphicx}
\usepackage{dcolumn}
\usepackage{bm}

\usepackage{mathptmx}
 \newtheorem{definition}{Definition}[section]

\newtheorem{example}{Example}[section]

\newcommand{\balpha}{\boldsymbol{\alpha}}

\newcommand{\bbeta}{\boldsymbol{\beta}}

\newcommand{\bSigma}{\boldsymbol{\Sigma} }
\newcommand{\bmu}{\boldsymbol{\mu} }

\newcommand{\bDelta}{\boldsymbol{\Delta} }
\newcommand{\bLambda}{\boldsymbol{\Lambda} }

\newcommand{\bgamma}{\boldsymbol{\gamma} }

\newcommand{\ba}{\boldsymbol{\alpha} }
\newcommand{\bepsilon}{\boldsymbol{\epsilon} }

\begin{document}
%

\title[Open Markov chains]{Open Markov chains: cumulant dynamics, fluctuations and correlations}

\author{R. Salgado-Garc\'{\i}a}
\email{raulsg@uaem.mx}
\affiliation{ 
Centro de Investigaci\'on en Ciencias - IICBA, Universidad Aut\'onoma del Estado de Morelos. Avenida Universidad 1001, colonia Chamilpa, C.P. 62209, Cuernavaca Morelos, Mexico.
}%

\date{\today}

\begin{abstract}
In this paper we propose a model for open Markov chains that can be interpreted as a system of non-interacting particles evolving according to the rules of a Markov chain. The number of particles in the system is not constant, because we allow the particles to arrive or leave the state space according to prescribed protocols. We describe this system by looking at the population of particles on every state by establishing the rules of time-evolution of the distribution of particles. We show that it is possible to describe the distribution of particles over the state space through the corresponding moment generating function. This description is given through the dynamics ruling the behavior of such a moment generating function and we prove that the system is able to attain the stationarity under some conditions. We also show that it is possible to describe the dynamics of the two first cumulants of the distribution of particles, which in some way is a simpler technique to obtain useful information of the open Markov chain for practical purposes. Finally we also study the behavior of the time-dependent correlation functions of the number of particles present in the system. We give some simple examples of open chains that either, can be fully described through the moment generating function or partially described through the exact solution of the cumulant dynamics.

\end{abstract}

\maketitle

\section{\label{sec:intro}Introduction}

Markov chains are discrete-time models for stochastic evolution, widely used to model systems in physics~\cite{lemons2002introduction,van1992stochastic}, chemistry~\cite{van1992stochastic}, biology~\cite{allen2010introduction,goel2016stochastic} among other disciplines~\cite{jackman2009bayesian,chib1996markov}. Roughly speaking, a Markov chain consists of a sequence of random variables $\{ X_t \in \mathcal{S} \, :\, t \in \mathbb{N}_0 \}$  taking values from a (finite or countable) set $\mathcal{S}$, called state space. The jump of the random variable from one state to another in one time step occurs with a prescribed probability, and the probabilities of all the possible jumps are collected in a matrix called \emph{Markov matrix} or \emph{stochastic matrix}. A natural way to interpret a Markov chain comes from physics; we can think of the random variable $X_t$ as the \emph{position at time} $t$ of given particle, and this particle moves on the discrete space $\mathcal{S}$. Thus, the stationary probability vector $\boldsymbol{\pi}$ (provided it exists) is interpreted from a point of view of \emph{ensembles}: if we have a collection of $N$ non-interacting particles moving according to the rules of the Markov chain, then, the stationary distribution of particles on $\mathcal{S}$ is $N\boldsymbol{\pi}$ if $N$ is large enough. This point of view clearly shows that a Markov chain  is a \emph{closed system}, since there is no inflow  of particles to $\mathcal{S}$ nor outflow of particles from $\mathcal{S}$. 

In this paper we shall be concerned with the case in which we allow the particles to arrive or leave the state space according to a prescribed \emph{protocol}. To be precise let us consider a state $j\in \mathcal{S}$. On one hand,  at every time step we allow a certain number of particles already present in the state $j$ to leave this state to the ``outside''. On the other hand we also allow a certain number of particles (from the ``outside'') to arrive at the state $j$. Both, the number of incoming particles and the number of outcoming particles, are modeled as random variables (or sequences of random variables) with a distribution given a priori. Our main goal in this paper is to describe the population of particles on the state space as well as its fluctuations. We are particularly interested in the behavior of both, the space correlations and time correlations for several possible scenarios for the incoming and outcoming protocols.  

At this point it is convenient to mention some works related to our model. For example, in a recent work~\cite{floriani2016flux}, Floriani \emph{et al} have considered the case in which some quantity (that they call ``mass'')  moves through the states of the chain according to the Markov transition probabilities. Particularly Floriani \emph{et al} focus their study to the case where the Markov chain has some absorbing states (in which the mass accumulates) and the mass is supplied by either, a  constant source or a periodic source. In contrast, in our model the number of incoming particles at every time step is not constant but a sequence of random variables not necessarily independent and identically distributed (i.i.d.). Moreover,  instead of modeling the outflow by absorbing states, we define a protocol of outcoming particles, which allows every particle to leave the chain with a state-dependent probability. Another work which is worth mentioning here is the one of Pollard and co-workers~\cite{baez2016compositional,pollard2016open,pollard2017open}.  They consider a class of open Markov process in which the state space is discrete and the time continuous. They assume that the incoming and outcoming fluxes are regulated  by means of a set of special states of $\mathcal{S}$ (which they call ``boundary''). One of the main differences of our approach with respect to the one proposed by Pollard is that they suppose that the elements of the boundary has a distribution prescribed a priori (which may be even time-dependent). In contrast, our approach considers every state as a source or sink of particles. In this sense  our model can be thought of as a Markov chain in contact with a \emph{reservoir of particles}. Thus, according to a prescribed protocol the particles go from the reservoir to the chain and viceversa, the particles go from the chain to the reservoir with a prescribed protocol. This way of modeling the source of particle is, in some way, similar to the grand canonical ensemble in thermodynamics, in which a system is in contact with a reservoir allowing particles to be interchanged.

\section{\label{sec:model}A model for open Markov chains}

The main idea behind our model for an ``open'' Markov chain  is that we allow the particles to \emph{arrive} and \emph{leave} the state space $\mathcal{S}$. We have mentioned that the particle can enter the state space according to a protocol which is modeled as a sequence of random variables. Such a sequence of random variables  determines the number of particles arriving at certain state every time step. On the other hand, the particle can leave the state space depending on which state they are. The most natural way to model this situation is by defining, for every state, a given probability with which the particle leaves such a state towards the reservoir. This probability must satisfy a compatibility condition, consisting in the fact that a particle in a given state has only two options i) jump to any other state in $\mathcal{S}$ or ii) jump to the reservoir. The sum of all these probabilities should be one, in order for the ``jump'' to be well defined. Notice that due to the compatibility condition, we have that the protocol of outcoming particles is completely determined by means of a non-negative matrix $\mathbf{Q}$ with a spectral radius strictly less than one. This is because the ``missing probability''  in $\mathbf{Q}$ (the necessary for $\mathbf{Q}$ to be a stochastic matrix) is interpreted as ``jump probabilities'' towards the outside.  

\begin{definition}[Open Markov chain]
\label{def:open_chain}
Let $\mathcal{S}$  be a finite set, whose cardinality is denoted by $\# \mathcal{S} = S$,  and let $\mathbf{Q} : \mathcal{S}\times  \mathcal{S} \to [0,1]\subset \mathbb{R}$ be an irreducible and aperiodic matrix with spectral radius strictly less than one. Let $\{ \mathbf{J}^{t}  : t \in \mathbb{N} \}$ be a sequence of random vectors taking values in $\mathbb{N}_0^{S}$. We say that $\big( \mathcal{S}, \mathbf{Q}, \{\mathbf{J}^t : t \in \mathbb{N}\} \big)$ is an open Markov chain with state space $\mathcal{S}$, jump matrix $\mathbf{Q}$ and incoming protocol $\{ \mathbf{J}^{t}  : t \in \mathbb{N} \}$. 
Now let $\{ \mathbf{N}^t : t \in \mathbb{R}\}$ be a sequence of random vectors taking values in $\mathbb{N}_0^{S}$. Such a sequence is defined as follows. Given the initial random vector $\mathbf{N}^0$ with a given distribution, we define $\mathbf{N}^t$ recursively as 
\begin{equation}
\label{eq:evolution}
\mathbf{N}^{t+1} := \mathbf{J}^t + \mathbf{R}^t,
\end{equation}
where $\mathbf{R}^t$ is a random vector taking values in $\mathbb{N}_0^{S}$, whose components are given by,
\begin{equation}
R^t_j := \sum_{i=1}^S B_{i,j}^t.
\end{equation}
The random variables $B_{i,j}^t$ are defined such that the $(S+1)$-dimensional vector $\mathbf{A}_i^t$, with components 
\[
(\mathbf{A}_{i}^t)_{j} =  \left\{ \begin{array} 
            {c@{\quad \mbox{ if } \quad}l} 
   B_{i,j}^t  &  1\leq j \leq  S  \\ 
   1- \sum_{j=1}^S B_{i,j}^t  &  j = S+1,  \\ 
             \end{array} \right. 
\]
has multinomial distribution, i.e., $ \mathbf{A}_i^t \sim \mbox{Multinomial} (\mathbf{z}_i,N^t_i) $. Also se assume that $ \mathbf{A}_i^t$ and $ \mathbf{A}_i^s$ are independent if $t\not= s$. Here $\{\mathbf{z}_i : 1\leq i \leq S\}$ is a set of probability vectors defined as
\[
(\mathbf{z}_i )_j =  \left\{ \begin{array} 
            {c@{\quad \mbox{ if } \quad}l} 
   q_{i,j}  &  1\leq j \leq  S  \\ 
   e_i &  j = S+1,  \\ 
             \end{array} \right. 
\]
where $q_{i,j}$ is the $(i,j)$th component of $\mathbf{Q}$ and $e_i$ is defined as,
\[
e_i :=1- \sum_{j=1}^S q_{i,j}.
\] 

We say that $\mathbf{N}^t$ is the \emph{distribution over the state space} at time $t$ with initial condition $\mathbf{N}^0$. 

\end{definition}

In the above definition we adopted the notation $\mathbb{N}_0$ to indicate the set $\{ 0,1,2,3,\dots \}$ to distinguish it from the set  $\mathbb{N} := \{ 1,2,3,\dots\}$. Notice  that  Eq.~(\ref{eq:evolution}) establishes the evolution of the number of particles. This equation states that the number of particles at time $t+1$ is the number of particles  having arrived as the state space (represented by $\mathbf{J}^t$) plus the number of particles having remained in the state space (which is represented by $\mathbf{R}^t$). Observe that the random vector $\mathbf{R}^t$ can be seen as a sum of independent random vectors $\mathbf{B}^t_i = (B_{i,j})_{j=1}^S $ representing the quantity of particles departing from state $i$ towards other states. Notice also that the random vector $\mathbf{A}_i^t$ is the ``enlarged'' version of the vector $\mathbf{B}_i^t$, because the vector $ \mathbf{A}_i^t$ represents the quantity of particles departing from state $i$ to other states and to the outside. It is clear that, due to conservation of particles during the process of redistribution, the vector $ \mathbf{A}_i^t$ should have multinomial distribution (its components actually sum up $N_i^t$ as we can see from the definition above). Hence we have that $\mathbf{N}^{t+1}$ depends on $ \mathbf{N}^{t}$ through the random variable $\mathbf{R}^t$ which is the responsible of the redistribution of particles among the internal states. 

From the above definition we can appreciate that the quantity of foremost interest is $\mathbf{N}^t$, the distribution of particles on the state space at time $t$. Our goal is then to provide a way for determining the probability function of $\mathbf{N}^t$ and to determine whether or not the process $\{\mathbf{N}^t:t\in \mathbb{N}_0\}$ is able to reach an stationary distribution. We will see that $\{\mathbf{N}^t : t\in \mathbb{N}_0\}$ is actually a Markov process and the main goal is to determine its properties. Before establishing a result on this direction let us give some examples of open Markov chains.

\begin{example}
\label{ex:one-vertex}
Let us consider the most simple case in which the system has only one state. In this case the matrix $\mathbf{Q}$ consists of a single number $q$ (a $1\times 1$ matrix), which should be non-negative and strictly less than one, i.e., $0\leq q <1$. The most simple case for the incoming protocol consists of a sequence of constant random variables all these taking the same value, which we call $J_0 \in \mathbb{N}_0$. This means that the number of particles arriving at every time-step is $J_0$. Every particle arriving at the unique available state has only two options jump to the outside or remain in its state. The probability of remaining in the state is $q$ and the probability of jumping to the outside is $1-q$. This simple example for open Markov chain can be illustrated as a graph with only one vertex and three edges as shown in 
figure~\ref{fig:one-vertex}. Notice that there are two special edges, one establishing the incoming protocol (labeled by $J_0$) and one defining the outcoming protocol (labeled by $1-q$).

\begin{figure}[ht]
\begin{center}
\scalebox{0.7}{\includegraphics{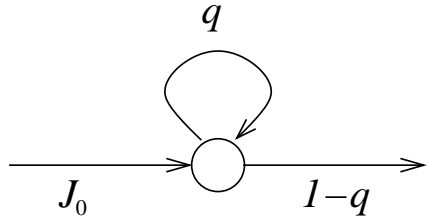}}
\end{center}
     \caption{
     One-vertex open chain. The circle represents the unique state available for the particles. The arrows in the graph stand for the ``jump'' rules allowed for the particles. Notice that the arrows that do not connect two states represent the incoming and outcoming protocols. 
     }
\label{fig:one-vertex}
\end{figure}
%

\end{example}

\begin{example}

A less trivial example is provided by giving a matrix larger than $1\times 1$. Let us consider for example the matrix given by 
\begin{equation}
 \mathbf{Q} = \left[
  \begin{array}{ccc}
  0 & 1/2 & 1/4 \\ 
1/4 & 1/4 & 0 \\ 
1/4 & 1/2 & 1/4
\end{array} \right]
\end{equation}
Notice that the above matrix is not stochastic, because some rows does not add to one, but less than one. The latter means that not all the states allow the particles to leave to the outside. As we can see the first row adds to $3/4$, which means that every particle on the state $1$ has a probability $1/4$ to  go out to the reservoir. The second row adds to $1/2$, this means that a given particle on the state $2$ has a probability $1/2$ to go out to the outside. Finally, the third row adds to one, meaning that a particle on the state $3$ can only jump to the other states $1$ or $2$ or remains in its current state, but it cannot leave the state space to go to the outside. 

Now assume that the incoming protocol $\{\mathbf{J}^t : t \in \mathbb{N}  \}$ is a set of i.i.d.~random vectors, i.e.~the protocol is time-independent. Then the number of incoming particles at every time-step can be considered as independent realizations of a single random vector, which we  denote by $\mathbf{J}$. To be more precise, we can chose in particular this random vector as $\mathbf{J}=(J_1,0,J_3)$, with $J_1$ and $J_2$ two independent random variables with Bernoulli distribution with parameters $p_1$ and $p_3$ respectively. If, for instance, the parameters of every Bernoulli distribution were given by $p_1=0.1$ and $p_3=0.6$, then this would mean that the average number of particles arriving at the vertex $1$ is lower than the average  number of particles arriving at the vertex $3$. Figure~\ref{fig:three-vertex} shows the graphical representation of this open Markov chain.

\begin{figure}[ht]
\begin{center}
\scalebox{0.5}{\includegraphics{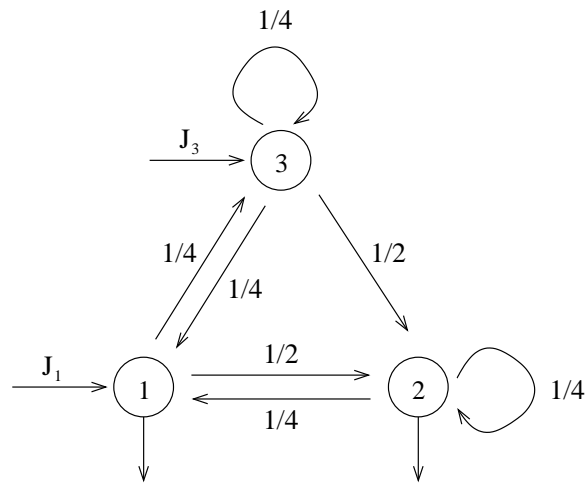}}
\end{center}
     \caption{
     An example of a three-states open Markov chain. 
     }
\label{fig:three-vertex}
\end{figure}
%

\end{example}

\section{\label{sec:evolution}The evolution of the particle distribution}


\subsection{Evolution of the moment generating function}


In this section we will establish the evolution equation for the process $\mathbf{N}^t$, which, as we have anticipated, is a Markov process. This fact can actually be appreciated in equation~(\ref{eq:evolution}), where we have indicated that the number of particles at time $t+1$ is uniquely determined by the number of incoming particles and the redistribution of particles that were present at time $t$. Clearly the last condition states the Markov property for the process  $\{\mathbf{N}^t : t\in \mathbb{N}_0\}$. In order to obtain the stochastic matrix governing the behavior of $\mathbf{N}^t$, let us define $p_t(\mathbf{n})$ as the probability vector associated to $\mathbf{N}^t$, i.e., 
\begin{equation}
p_t(\mathbf{n}) := \mathbb{P}( \mathbf{N}^t = \mathbf{n}).
\end{equation}
We will refer to $p_t(\mathbf{n})$ as the \emph{distribution of particles over the state space} or, to simplify,  \emph{distribution over the state space}
Let us consider the probability vector at time $t+1$. Notice that, due the fact that $\mathbf{N}^{t+1}$ depends only on $\mathbf{N}^t$, which is a consequence of the relation~(\ref{eq:evolution}), as we have mentioned above, it is clear that
\begin{eqnarray}
p_{t+1} (\mathbf{n}) &=&
 \mathbb{P} (\mathbf{N}^{t+1} = \mathbf{n}) 
\nonumber
\\
&=& \sum_{\mathbf{k} \in \mathbb{N}_0^S} \mathbb{P} (\mathbf{N}^{t} = \mathbf{k}) 
\mathbb{P} (\mathbf{N}^{t+1} = \mathbf{n}| \mathbf{N}^{t} = \mathbf{k} ). 
\label{eq:1st-evolution}
\end{eqnarray}
Observe that the last expression is a consequence of the Markov property of the process which was assumed in equation~(\ref{eq:evolution}). Notice that equation~(\ref{eq:1st-evolution}) establishes the evolution equation for $p_t$, which can be written as
\begin{eqnarray}
p_{t+1} (\mathbf{n})  &=& \sum_{\mathbf{k} \in \mathbb{N}_0^S} p_t(\mathbf{k}) 
\mathbb{P} (\mathbf{N}^{t+1} = \mathbf{n}| \mathbf{N}^{t} = \mathbf{k} ). 
\end{eqnarray}
Now let us call $K(\mathbf{k},\mathbf{n})$ the conditional probability appearing in the above equation, i.e., 
\begin{equation}
K(\mathbf{k},\mathbf{n}) := \mathbb{P} (\mathbf{N}^{t+1} = \mathbf{n}| \mathbf{N}^{t} = \mathbf{k} ), 
\end{equation}
and observe that this quantity can be rewritten as follows
\begin{equation}
K(\mathbf{k},\mathbf{n}) = \mathbb{P} (\mathbf{J}^{t} + \mathbf{R}^{t}  = \mathbf{n}).
\end{equation}
The dependence on $\mathbf{k}$ in the above expression is implicit in the random vector $\mathbf{R}^t$, since the redistribution of particles depends on the number of particles on every state at time $t$, which is indeed given  by $\mathbf{k}$. Thus, the function $K : \mathbb{N}_0^S \times \mathbb{N}_0^S \to [0,1] \subset\mathbb{R} $ can be thought  as the \emph{stochastic matrix} corresponding to the process $\{\mathbf{N}^t : t \in \mathbb{N}_0\}$, since the evolution equation for the probability vector $p_t(\mathbf{n})$ is given in terms of $K$ as follows,
\begin{equation}
\label{eq:evolution_pK}
p_{t+1} (\mathbf{n})  = \sum_{\mathbf{k} \in \mathbb{N}_0^S} p_t(\mathbf{k}) 
K (\mathbf{k}, \mathbf{n} ). 
\end{equation}

In order to solve equation~(\ref{eq:evolution_pK}) for $p_t (\mathbf{n})$ it is necessary to make some assumptions on the nature of the random vectors $\mathbf{J}^t$ and $\mathbf{R}^t$ for all $t$.  First of all, it is natural to assume that $\mathbf{J}^t$ and $\mathbf{R}^t$ are independent. This assumption actually means that the number of particles incoming to the state space does not have any influence on the redistribution of particles already present in the chain. This implies that the joint probability for the random vectors $\mathbf{J}^t$ and $\mathbf{R}^t$ can be factorized as the product of its corresponding probability vectors, i.e., 
\[
\mathbb{P} (\mathbf{J}^t = \mathbf{j} ; \mathbf{R}^t = \mathbf{r}) = \mathbb{P} (\mathbf{J}^t = \mathbf{j}) \mathbb{P} ( \mathbf{R}^t = \mathbf{r}).
\] 
The above equality allows us to express the kernel $K(\mathbf{k},\mathbf{n})$ as,
\begin{eqnarray}
K(\mathbf{k},\mathbf{n}) &:=& \mathbb{P} (\mathbf{J}^{t} + \mathbf{R}^{t}  = \mathbf{n}) 
\nonumber
\\
&=& \sum_{\mathbf{j} +\mathbf{r} = \mathbf{n}}  \mathbb{P} (\mathbf{J}^t = \mathbf{j}) \mathbb{P} ( \mathbf{R}^t = \mathbf{r}).
\label{eq:K-sum}
\end{eqnarray}

Next, we will solve for $p_t$ by using the well-known technique of the \emph{moment generating function} (m.g.f.). To this end let us introduce some notation. Let $\mathcal{G}_t : \mathbb{R}^S \times \mathbb{R}^S \to \mathbb{R} $ be the m.g.f. of $\mathbf{N}^t$, which is defined as
\begin{equation}
\mathcal{G}_t (\boldsymbol{\alpha}) := \mathbb{E} \left[  e^{\mathbf{N}^t} \balpha^T\right] 
= \sum_{\mathbf{n}\in\mathbb{N}_0^S } p_t(\mathbf{n}) e^{\mathbf{n}\balpha^T  }.
\end{equation}  
At this point it is important to describe our convention for vectors in $\mathbb{R}^S$. First of all we should emphasize that we interpret the vectors  $\mathbf{n}$, $\balpha $,  $\mathbf{N}^t$, etc., as \emph{row} vectors (i.e. matrices of size $1\times S$). Thus, the superscript $T$ means, as usual, matrix transposition implying that the vector  $\boldsymbol{\alpha}^T$ is a column vector (a matrix of size $S\times 1$). Thus, within this convention, the product $\mathbf{n}\balpha^T$ should be understood in the sense of the usual matrix product, which in this case results in a single number.

Analogously, we also define the moment generating functions for $\mathbf{J}^t$ and $\mathbf{R}^t$ as follows,
\begin{eqnarray}
\mathcal{F}_t(\balpha) &:=&   \mathbb{E} \left[  e^{\mathbf{J}^t\balpha^T}\right] 
= \sum_{\mathbf{j}\in\mathbb{N}_0^S }   \mathbb{P} (\mathbf{J}^t = \mathbf{j})  e^{ \mathbf{j}\balpha^T },
\\
\mathcal{H} (\ba) &:=& \mathbb{E} \left[  e^{\mathbf{R}^t\ba^T}\right] 
= \sum_{\mathbf{r}\in\mathbb{N}_0^S }  \mathbb{P} (\mathbf{R}^t = \mathbf{r})  e^{\mathbf{r}\balpha^T }.
\end{eqnarray}
Notice that we omitted the superscript $t$ in the m.g.f.~for $\mathbf{R}^t$, since two random vectors, say for example  $\mathbf{R}^t$ and  $\mathbf{R}^{s}$, are independent and have the same distribution, and consequently, share the same moment generating function.

Next, our objective is to provide a recurrence relation for $\mathcal{G}_t$ using the evolution equation~(\ref{eq:evolution_pK}). Thus, let us consider the m.g.f.~for $\mathbf{N}^{t+1}$, and note that
\begin{eqnarray}
\mathcal{G}_{t+1} (\ba)&=& \sum_{\mathbf{n}\in\mathbb{N}_0^S } p_{t+1}(\mathbf{n}) e^{ \mathbf{n} \balpha^T }
\nonumber
\\
&=& \sum_{\mathbf{n}\in\mathbb{N}_0^S }  \sum_{\mathbf{k} \in \mathbb{N}_0^S} p_t(\mathbf{k}) K (\mathbf{k}, \mathbf{n} ) e^{\mathbf{n} \balpha^T}
\nonumber
\\
&=& \sum_{\mathbf{n}\in\mathbb{N}_0^S }  \sum_{\mathbf{k} \in \mathbb{N}_0^S}  \sum_{\mathbf{j} +\mathbf{r} = \mathbf{n}} 
p_t(\mathbf{k})  \mathbb{P} (\mathbf{J}^t = \mathbf{j}) \mathbb{P} ( \mathbf{R}^t = \mathbf{r})  e^{\mathbf{n} \balpha^T},
\end{eqnarray}
where we have used the form of $K$ given in equation~(\ref{eq:K-sum}). We can appreciate that the summation over $\mathbf{n}$ together with the summation over the restriction $ \mathbf{j} +\mathbf{r} = \mathbf{n}$ results in a double sum over the ``indices'' $\mathbf{j}$ and $ \mathbf{r} $ without restrictions, i.e., we obtain two sums over independent indices. This observation allows us to write
\begin{eqnarray}
\mathcal{G}_{t+1} (\ba)  &=&   \sum_{\mathbf{k} \in \mathbb{N}_0^S}  \sum_{\mathbf{j}\in\mathbb{N}_0^S  }  \sum_{\mathbf{r}\in\mathbb{N}_0^S } 
p_t(\mathbf{k})  \mathbb{P} (\mathbf{J}^t = \mathbf{j}) \mathbb{P} ( \mathbf{R}^t = \mathbf{r})  e^{ (\mathbf{j} + \mathbf{r}) \balpha^T }.
\end{eqnarray}
In the above equation we can observe that the summation over $\mathbf{j}$ and $\mathbf{r}$ results in the m.g.f.~for $\mathbf{J}^t$ and $\mathbf{R}^t$ respectively. This implies that
\begin{eqnarray}
\mathcal{G}_{t+1} (\ba)  &=&   \sum_{\mathbf{k} \in \mathbb{N}_0^S} 
p_t(\mathbf{k}) \,   \mathcal{F}_t(\ba)\, \mathcal{H}(\ba) .
\label{eq:Gt+1-sum}
\end{eqnarray}

In Appendix~\ref{ape:1} we show that $\mathcal{H}(\ba)$ can be written as
\begin{equation}
\label{eq:gen-H-k}
\mathcal{H}(\ba) = e^{\mathbf{k} \mathbf{H}^T(\ba) },
\end{equation}
where the function $\mathbf{H} : \mathbb{R}^S \to \mathbb{R}^S$ is defined as follows. If $H_i (\ba) =\left( \mathbf{H} (\ba) \right)_i $ denotes the $i$th component of $\mathbf{H}$ we have that
\begin{equation}
H_i(\ba) :=  \log \left( e_i +  \sum_{j=1}^S q_{i,j} e^{\alpha_j}  \right).
\end{equation}
Observe that the equation~(\ref{eq:gen-H-k}) shows explicitly the dependence on $\mathbf{k}$ of the m.g.f. for $\mathbf{R}^t$. Then, if we substitute the relation~(\ref{eq:gen-H-k}) into~(\ref{eq:Gt+1-sum}), we obtain, 
\begin{eqnarray}
\mathcal{G}_{t+1} (\ba)  &=&   \sum_{\mathbf{k} \in \mathbb{N}_0^S}  p_t(\mathbf{k}) \,   \mathcal{F}_t(\ba) e^{\mathbf{k} \mathbf{H}^T(\ba)}.
\end{eqnarray}
We can easily see that the summation over $\mathbf{k}$ results in the m.g.f.~for $\mathbf{N}^{t}$.  Thus,
\begin{eqnarray}
\mathcal{G}_{t+1} (\ba)  &=&  \mathcal{F}_t(\ba) \, \mathcal{G}_{t} \left( \mathbf{H}(\ba) \right). 
\label{eq:Gt+1}
\end{eqnarray}

Equation~(\ref{eq:Gt+1}) is a recurrence relation governing the time-dependence of the m.g.f.~for $\mathbf{N}^t$. This equation can be formally solved to obtain 
\begin{eqnarray}
\mathcal{G}_{t} (\ba)  &=&  \mathcal{G}_{0} \left(\mathbf{H}^{(t+1)}(\ba) \right) \, \prod_{r=0}^{t-1} \mathcal{F}_{t-r}\left( \mathbf{H}^{(r)}(\ba) \right),
\label{eq:Gt_formal-sol}
\end{eqnarray}
where $\mathcal{G}_{0}$ stand for the m.g.f for $\mathbf{N}^0$  (the initial distribution on the state space). We should emphasize that the superscript notation $\mathbf{H}^{(r)}$ stands for the $r$th iteration of the function $\mathbf{H}$, i.e., $\mathbf{H}^{(r)} := \mathbf{H} \circ \mathbf{H} \circ \dots \mathbf{H} $, $r$ times.

Now, let us assume that the process $\{ \mathbf{J}^t : t\in \mathbb{N} \}$ is a sequence of identically distributed random vectors (not necessarily independent). In this case, the m.g.f. for $\mathbf{J}^t$ is the same for all $t$, consequently the formal solution for $\mathcal{G}_t$ can be expressed as,
\begin{eqnarray}
\mathcal{G}_{t} (\ba)  &=&  \mathcal{G}_{0} \left(\mathbf{H}^{(t+1)}(\ba) \right) \, \prod_{r=0}^{t-1} \mathcal{F}\left( \mathbf{H}^{(r)}(\ba) \right).
\label{eq:Gt_formal-sol-iid-J}
\end{eqnarray}
This result, together with the fact that $\mathbf{H}^{(r)} (\balpha) \to \mathbf{0}$ as $r \to \infty$ in an open neighborhood around $\balpha = \mathbf{0}$ (see Appendix~\ref{ape:1} for a proof), implies that, for the case in which the random vectors $\mathbf{J}^t$  are identically distributed for all $t$,  the process $\{ \mathbf{N}^t : t\in \mathbb{N}_0\}$ admits an stationary solution. This is because the m.g.f. $\mathcal{G}_t$ attains a limit when $t\to \infty$. Such a limit can be written as,
\begin{eqnarray}
\mathcal{G}_{\mathrm{stat}} (\ba)  &=&   \prod_{r=0}^{\infty} \mathcal{F}\left( \mathbf{H}^{(r)}(\ba) \right).
\label{eq:Gt_steady}
\end{eqnarray}
whenever the infinite product exist. If this is the case, the m.g.f. $\mathcal{G}_{\mathrm{stat}} (\ba) $ is additionally a solution for the evolution equation~(\ref{eq:Gt+1}). This means that $\mathcal{G}_{\mathrm{stat}} (\ba) $ corresponds to a m.g.f. of a distribution $p_{\mathrm{stat}} (\mathbf{n})$ over the state space, which is invariant under the dynamics~(\ref{eq:evolution_pK}).

\begin{example}
\label{ex:one-vertex-2}
Let us consider the open chain consisting of only one vertex given in example~\ref{ex:one-vertex}. In this case we will assume that the incoming protocol $\{J^t : t\in\mathbb{N}\}$ consists of a sequence of i.i.d. random variables having Bernoulli distribution with parameter $p$. Since all $J^t$ have identical distribution, then we have only one m.g.f. $\mathcal{F}$ characterizing them. This function is given by,
\[
\mathcal{F} (\alpha) =  1-p + p e^{\alpha}.
\] 
On the other hand, due to the fact that  $\mathbf{Q}$ is a $1\times 1$ matrix, we have that $\mathbf{H}$  is a real-valued function depending on one variable, $\alpha$, given by,
\[
H (\alpha) =  \log(1-q + q e^{\alpha}).
\] 
where we denoted by $q$ the unique element of $\mathbf{Q}$. In this case, the $r$th iteration of $H$ can be exactly calculated. A straightforward calculation shows that,
\[
H^{(r)}(\alpha) = \log \left( 1-q^r + q^r e^{\alpha}\right) .
\]
Notice that $H^{(r)}(\alpha) \to 0$ as $r\to \infty$, as we have anticipated above. Hence, in this case, the stationary m.g.f.~is given by
\begin{eqnarray}
\mathcal{G}_{\mathrm{stat}} (\alpha)  &=&   \prod_{r=0}^{\infty} \mathcal{F}\left( {H}^{(r)}(\alpha) \right) =  \prod_{r=0}^{\infty}  \left( 1-p + p \left( 1-q^r + q^r e^{\alpha}\right)\right)
\nonumber
\\
 &=&  \prod_{r=0}^{\infty}  \left( 1 - p q^r + p q^r e^{\alpha} \right).
\end{eqnarray}
The above result shows that the stationary distribution can be interpreted as the convolution of an infinite sequence of Bernoulli distributions with parameters $p q^r$. This particularly means that the random variable $N^t$, when it has reached the stationarity can be written as a sum of Bernoulli random variables $X_r$ (with parameter $pq^r$),
\[
N^t = \sum_{r=0}^\infty X_r.
\]
The above result allows us to compute, for example, the expected value $\mathbb{E}[N^t] $ and the variance,
\begin{eqnarray}
\mathbb{E}[N^t] &=& \sum_{r=0}^\infty \mathbb{E}[X_r] = \sum_{r=0}^\infty pq^r = \frac{p}{1-q},
\label{eq:1v_mean}
\\
\mathrm{Var}(N^t) &=& \sum_{r=0}^\infty \mathrm{Var}(X_r)  = \sum_{r=0}^\infty  pq^r (1-pq^r) = \frac{p}{1-q} + \frac{p^2}{1-q^2}.
\label{eq:1v_var}
\end{eqnarray}

We should observe that the fact that the expected value $\mathbb{E}[N^t] = p/(1-q)$ is finite implies that an equilibration is attained between the number of incoming  particles  and the number of outcoming particles. It is clear that the mean number of arriving particles per unit time is $p$, while, the number of leaving particles per unit time is $ (1-q) \mathbb{E}[N^t] $. Once $N^t$ has  reached the stationarity, we have an equality between these quantities, giving the result stated above. Although we have obtained the mean number of particles by using the argument of equilibration, the same line of reasoning cannot be applied to the variance. We have thus provided a way to compute the fluctuations in the number of particles that are present in the vertex.

\end{example}


\subsection{Cumulant dynamics} 


Up to now we have obtained two main results, a recurrence relation of the time-dependent m.g.f. of $\mathbf{N}^t$ and  a formal expression for the m.g.f. of $\mathbf{N}^t$ when the system has attained the stationarity. Now we will focus in two quantities which are of special interest, namely, the two first cumulants for $\mathbf{N}^t$, and how do these quantities evolve in time. Let us start by obtaining the first cumulant of $\mathbf{N}^t$. Notice that the first cumulant (which coincides with the first moment) can be obtained by taking the first derivative of the m.g.f.~$ \mathcal{G}_t (\alpha)$. Let us denote by $\bmu_t$ the first moment of $\mathbf{N}^t$, i.e., 
\begin{equation}
\bmu_t := \mathbb{E} [\mathbf{N}^t],
\end{equation}
to which we will refer to as the \emph{mean distribution of particles over the state space} at time $t$, or simply,  the \emph{mean distribution} at time $t$. Now let us notice that, 
\begin{equation}
(\bmu_t)_i = \frac{\partial \mathcal{G}_t(\balpha)}{\partial \alpha_i}\bigg|_{\balpha = \mathbf{0}},
\end{equation}
where $(\bmu_t)_i$ is the $i$th component of $\bmu_t$  and $\alpha_i$ is the $i$th component of $\balpha$. Observe that the $i$th component of the mean distribution at time $t +1$ can be obtained from equation~(\ref{eq:Gt+1}), giving
\begin{eqnarray}
(\bmu_{t+1})_i &=& \frac{\partial \mathcal{G}_{t+1}(\balpha)}{\partial \alpha_i}\bigg|_{\balpha = \mathbf{0}}
\nonumber
\\
&=&
\frac{\partial \mathcal{F}_t (\balpha)}{\partial \alpha_i}  \mathcal{G}_{t} \left( \mathbf{H}(\ba) \right)\bigg|_{\balpha = \mathbf{0}} + 
\mathcal{F}_t (\balpha) \sum_{k=1}^S \frac{\partial  \mathcal{G}_{t} \left( \mathbf{H} \right) }{ \partial H_k} \frac{ \partial H_k}{\partial \alpha_i }\bigg|_{\balpha = \mathbf{0}}.
\end{eqnarray}
Notice that $\mathbf{H}(\balpha = \mathbf{0} ) = \mathbf{0}$ and that any moment generating function evaluated at $\mathbf{0}$ is one. Hence we have, 
\begin{eqnarray}
\label{eq:mut-1}
(\bmu_{t+1})_i &=&
(\bepsilon_t)_i   + \sum_{k=1}^S (\bmu_t)_k q_{k,i},
\end{eqnarray}
where $\bepsilon_t$ stands for the expected value of $\mathbf{J}^t$, which can be obtained by means of the first derivative of $\mathcal{F}_t$, i.e.,
\[
\frac{\partial \mathcal{F}_t (\balpha)}{\partial \alpha_i} \bigg|_{\balpha = \mathbf{0}} = (\bepsilon_t)_i.
\]
In equation~(\ref{eq:mut-1}) made use of  the fact that
\[
 \frac{ \partial H_k}{\partial \alpha_i }\bigg|_{\balpha = \mathbf{0}} = (\mathbf{Q})_{k,i} = q_{k,i},
\]
which is proved in~\ref{ape:1}. Thus, it is clear that equation~(\ref{eq:mut-1}) can be written as,
\begin{equation}
\bmu_{t+1} = \bepsilon_t   + \bmu_t \mathbf{Q}.
\label{eq:mut+1}
\end{equation}
The above expression states the dynamics for the evolution of the mean distribution $\bmu_t$ in time. This evolution has two components, one involving the internal dynamics (which is given by the term $\bmu_t \mathbf{Q}$ giving the internal redistribution of particles) and other one involving the external dynamics (which is given by the time-dependent mean  number of incoming particles).

Now let us explore the behavior of the second cumulant of $\mathbf{N}^t$. The second cumulant corresponds to the variance matrix $\mathrm{Var}(\mathbf{N}^t)$ which we will denote by $\bSigma_t$. Notice that this matrix has entries given by
\begin{equation}
(\bSigma_t)_{i,j} := \left( \mathrm{Var}(\mathbf{N}^t) \right)_{i,j} = \mathbb{E}[N^t_i N^t_j] - (\bmu_t)_i (\bmu_t)_j. 
\end{equation}
Next we will use the dynamics of the m.g.f.~of $\mathbf{N}^t$ to obtain a recurrence for the expected value $ \mathbb{E}[N^t_i N^t_j] $. It is clear that
\begin{equation}
 \mathbb{E}[N^{t+1}_i N^{t+1}_j]  = \frac{\partial^2 \mathcal{G}_{t+1} (\balpha)}{\partial \alpha_i \partial \alpha_j}  \bigg|_{\balpha = \mathbf{0}}.
\end{equation}
The above expression together with the evolution equation~(\ref{eq:Gt+1}) leads to
\begin{eqnarray}
\mathbb{E}[N^{t+1}_i N^{t+1}_j]  &=&   \frac{\partial^2}{\partial \alpha_i \partial \alpha_j} \bigg( \mathcal{F}_t(\ba) \, \mathcal{G}_{t} \left( \mathbf{H}(\ba) \right) \bigg) \bigg|_{\balpha = \mathbf{0}}
\nonumber
\\
&=&
\frac{\partial}{\partial \alpha_i } 
\bigg( 
\frac{\partial \mathcal{F}_t (\balpha)}{\partial \alpha_j}  \mathcal{G}_{t} \left( \mathbf{H}(\ba) \right)  + 
\mathcal{F}_t (\balpha) \sum_{k=1}^S \frac{\partial  \mathcal{G}_{t} \left( \mathbf{H} \right) }{ \partial H_k} \frac{ \partial H_k}{\partial \alpha_j }   \bigg) \bigg|_{\balpha = \mathbf{0}}
\nonumber
\\
&=&
\bigg( \frac{\partial^2 \mathcal{F}_t (\balpha)}{\partial \alpha_i \partial \alpha_j}  \mathcal{G}_{t} \left( \mathbf{H}(\ba) \right)  +
\frac{\partial \mathcal{F}_t (\balpha)}{\partial \alpha_j} \sum_{k=1}^S \frac{\partial  \mathcal{G}_{t} \left( \mathbf{H} \right) }{ \partial H_k} \frac{ \partial H_k}{\partial \alpha_i } 
\nonumber
\\
&+& 
\frac{\partial \mathcal{F}_t (\balpha)}{\partial \alpha_i} \sum_{k=1}^S \frac{\partial  \mathcal{G}_{t} \left( \mathbf{H} \right) }{ \partial H_k} \frac{ \partial H_k}{\partial \alpha_j }  
+ 
 \mathcal{F}_t (\balpha) 
\sum_{k=1}^S \sum_{l=1}^S \frac{\partial^2  \mathcal{G}_{t} \left( \mathbf{H} \right) }{\partial H_l \partial H_k} \frac{ \partial H_k}{\partial \alpha_i } \frac{ \partial H_l}{\partial \alpha_j }  
\nonumber
\\
&+&
 \mathcal{F}_t (\balpha)\sum_{k=1}^S
 \frac{\partial  \mathcal{G}_{t} \left( \mathbf{H} \right) }{ \partial H_k} \frac{ \partial^2 H_k}{\partial \alpha_i  \partial \alpha_j } 
 \bigg) \bigg|_{\balpha = \mathbf{0}}.
 \label{eq:varNt-1}
\end{eqnarray}
At this point it is necessary to introduce some notations. First let us denote by $\bDelta$ the variance matrix of the incoming flux, i.e.,
\begin{equation}
\bDelta_t := \mathrm{Var} (\mathbf{J}^t).
\end{equation}
The above quantity can be obtained through the second derivative of the m.g.f. $\mathcal{F}_t(\balpha)$ as follows,
\begin{equation}
(\bDelta_t)_{i,j} =  \mathbb{E}[J^t_i J^t_j ] -  \mathbb{E}[J^t_i ] \mathbb{E}[J^t_j] = \frac{\partial^2 \mathcal{F}_t (\balpha)}{\partial \alpha_i \partial \alpha_j} - (\bepsilon_t)_i (\bepsilon_t)_j
\end{equation}
Thus, performing the evaluation of equation~(\ref{eq:varNt-1}) at $\balpha = \mathbf{0}$ we obtain,
\begin{eqnarray}
\mathbb{E}[N^{t+1}_i N^{t+1}_j]  &=&  (\bDelta_t)_{i,j} + (\bepsilon_t)_i(\bepsilon_t)_j + (\bepsilon_t)_j \sum_{k=1}^S (\bmu_t)_k q_{k,i} 
+(\bepsilon_t)_i \sum_{k=1}^S (\bmu_t)_k q_{k,j}
\nonumber
\\
&+&
\sum_{k=1}^S \sum_{l=1}^S \mathbb{E}[N^{t}_k N^{t}_l] q_{k,i} q_{l,j} + 
\sum_{k=1}^S (\bmu_t)_k \left( q_{k,i}\delta_{i,j}-q_{k,i}q_{k,j} \right),
\label{eq:varNt-2}
\end{eqnarray}
where we used the fact that (see Appendix~\ref{ape:1}),
\begin{equation}
\frac{\partial^2 H_k (\balpha)}{ \partial \alpha_i \partial \alpha_j } \bigg|_{\balpha = \mathbf{0}} =
q_{k,i}\delta_{i,j}-q_{k,i}q_{k,j}.
\end{equation}
Let us simplify equation~(\ref{eq:varNt-2}) by noticing that the summations can be written as matrix products,
\begin{eqnarray}
\mathbb{E}[N^{t+1}_i N^{t+1}_j]  &=&  (\bDelta_t)_{i,j} + (\bepsilon_t)_i(\bepsilon_t)_j 
+ (\bepsilon_t)_j (\bmu_t \mathbf{Q})_i + (\bepsilon_t)_i  (\bmu_t \mathbf{Q})_j + (\bLambda_t)_{i,j}
\nonumber
\\
&+&
\sum_{k=1}^S \sum_{l=1}^S \mathbb{E}[N^{t}_k N^{t}_l] q_{k,i} q_{l,j},
\label{eq:varNt-3}
\end{eqnarray}
where we defined the matrix $\bLambda_t$ as
\begin{equation}
\label{eq:Lambda_t}
(\bLambda_t)_{i,j} := \sum_{k=1}^S (\bmu_t)_k \left( q_{k,i}\delta_{i,j}-q_{k,i}q_{k,j}  \right).
\end{equation}
Equation~(\ref{eq:varNt-3}) allows us to obtain the variance matrix $\mathrm{Var}(\mathbf{N}^{t+1})$,
\begin{eqnarray}
\mathrm{Var} (\mathbf{N}^{t+1}) &=&  \mathbb{E}[N^{t+1}_i N^{t+1}_j] - (\bmu_{t+1})_i (\bmu_{t+1})_j 
\nonumber
\\
&=&
 (\bDelta_t)_{i,j} + (\bLambda_t)_{i,j}  + (\bepsilon_t)_i(\bepsilon_t)_j  + 
 (\bepsilon_t)_j (\bmu_t \mathbf{Q})_i + (\bepsilon_t)_i  (\bmu_t \mathbf{Q})_j 
\nonumber
\\
&+&
 \sum_{k=1}^S \sum_{l=1}^S \left( \mathbb{E}[N^{t}_k N^{t}_l] -  (\bmu_t)_k (\bmu_t)_l \right)  q_{k,i} q_{l,j} +  \sum_{k=1}^S \sum_{l=1}^S   (\bmu_t)_k (\bmu_t)_l   q_{k,i} q_{l,j} 
 \nonumber
 \\
&-& (\bmu_{t+1})_i (\bmu_{t+1})_j 
\nonumber
\\
&=&
 (\bDelta_t)_{i,j} + (\bLambda_t)_{i,j}  + (\bepsilon_t)_i(\bepsilon_t)_j  + 
 (\bepsilon_t)_j (\bmu_t \mathbf{Q})_i + (\bepsilon_t)_i  (\bmu_t \mathbf{Q})_j 
\nonumber
\\
&+&
 \sum_{k=1}^S \sum_{l=1}^S q_{k,i} \left( \mathrm{Var} (\mathbf{N}^t) \right)_{k,l}   q_{l,j} + 
  (\bmu_t \mathbf{Q})_i(\bmu_t \mathbf{Q})_j -(\bmu_{t+1})_i (\bmu_{t+1})_j .
\label{eq:varNt-4}
\end{eqnarray}
Rearranging terms in the above expression and denoting by $\bSigma_t$ the matrix variance $\mathrm{Var}(\mathbf{N}^t)$, we obtain 
\begin{eqnarray}
(\bSigma_{t+1})_{i,j}
&=&
 (\bDelta_t+\bLambda_t)_{i,j}  + (\bepsilon_t)_i(\bepsilon_t)_j  + 
 (\bepsilon_t)_j (\bmu_t \mathbf{Q})_i + (\bepsilon_t)_i  (\bmu_t \mathbf{Q})_j  +   (\bmu_t \mathbf{Q})_i(\bmu_t \mathbf{Q})_j
\nonumber
\\
&+&
 \sum_{k=1}^S \sum_{l=1}^S q_{k,i} (\bSigma_{t})_{k,l}   q_{l,j} - (\bmu_{t+1})_i (\bmu_{t+1})_j 
 \nonumber
 \\
 &=&
 (\bDelta_t+\bLambda_t)_{i,j} + \left( \bepsilon_t + \bmu_t \mathbf{Q} \right)_i \left( \bepsilon_t + \bmu_t \mathbf{Q} \right)_j
+ \left( \mathbf{Q}^T \bSigma_t \mathbf{Q} \right)_{i,j} - (\bmu_{t+1})_i (\bmu_{t+1})_j.
\nonumber
\\
\label{eq:varNt-5}
\end{eqnarray}
Finally, observing that $\bepsilon_t + \bmu_t \mathbf{Q}  = \bmu_{t+1}$, it is clear that the variance matrix satisfy the evolution equation,
\begin{equation}
\bSigma_{t+1} =  \bDelta_t+\bLambda_t + \mathbf{Q}^T \bSigma_t \mathbf{Q}.
\label{eq:varNt+1}
\end{equation}

Equations~(\ref{eq:mut+1}) and~(\ref{eq:varNt+1}) govern the dynamics of the first and second cumulants and are valid even when the incoming flux has a time-dependent distribution. In the case where the protocol of incoming particles $\{ \mathbf{J}^t : t \in\mathbb{N}\}$ is a stationary process (for which the two first cumulants are time-independent) we have that the system can reach the stationarity.  Particularly we have that the dynamics equations~(\ref{eq:mut+1}) and~(\ref{eq:varNt+1}) have stationary solutions (proved in~\ref{ape:1}) given by,
\begin{eqnarray}
\label{eq:mu_stat} 
\overline{\bmu} &=& \bepsilon (\mathbf{1} - \mathbf{Q})^{-1}
\\
\label{eq:sigma_stat}
\overline{\bSigma} &=& \sum_{k=0}^\infty (\mathbf{Q}^T)^k  (\bDelta+\overline{\bLambda}  )\mathbf{Q}^k,
\end{eqnarray}
where $\bepsilon$ and $\bDelta $ are the mean vector and the variance matrix of the stationary process $\{ \mathbf{J}^t : t \in\mathbb{N}\}$ and  $ \overline{\bmu}$ and $ \overline{\bSigma}$ denote the mean distribution $\mathbb{E}[\mathbf{N}^t]$ and the variance matrix $\mathrm{Var} (\mathbf{N}^t)$ when the process  $\{ \mathbf{N}^t : t \in\mathbb{N}_0\}$ has reached the stationarity.  We also defined  $\overline{\bLambda} $ as 
\begin{equation}
\label{eq:barLambda}
( \overline{\bLambda} )_{i,j} = \sum_{k=0}^S (\overline{\bmu})_k \left( q_{k,i}\delta_{i,j}-q_{k,i}q_{k,j}  \right) . 
\end{equation}

\begin{example}
\label{ex:3v_symmetric}

Let us consider a three states open chain  with jump matrix given by
\begin{equation}
 \mathbf{Q} = \left[
  \begin{array}{ccc}
  0 & q & q \\ 
q & 0 & q \\ 
q & q & 0
\end{array} \right]
\end{equation}
with $p$ a parameter restricted to take values in the interval $0 < q < 1/2$. Let us also assume that the protocol of incoming particles  $\{ \mathbf{J}^t = (J_1^t,J_2^t,J_3^t) : t \in\mathbb{N}\}$  is a stationary process whose joint probability distribution at time $t$ can be written as  
\[
\mathbb{P}(J_1^t = j_1;J_2^t = j_2;J_3^t = j_3) =: f_{1,2}(j_1, j_2) f_3(j_3).
\] 
Particularly, we chose $f_{1,2}$ and $f_3$ as
\begin{eqnarray}
f_{1,2} (j_1, j_2) &=&  \left\{ \begin{array} 
            {r@{\quad \mbox{ if } \quad}l} 
   p/2  &  (j_1,j_2) = (1,1)  \\ 
   p/2  &  (j_1,j_2) = (0,0)  \\ 
   (1-p)/2  &  (j_1,j_2) = (1,0)  \\ 
   (1-p)/2  &  (j_1,j_2) = (0,1)  \\ 
   0  & \mathrm{otherwise}  \\ 
             \end{array} \right. 
\\    
f_{3} (j_3) &=&  \left\{ \begin{array} 
            {r@{\quad \mbox{ if } \quad}l} 
   1/2  &  j_3=0  \\ 
   1/2  &  j_3=1  \\ 
   0  & \mathrm{otherwise}  \\ 
             \end{array} \right. 
\end{eqnarray}
Notice that each random variable $J_i^t$ can only take the values $0$ or $1$. Moreover, the above choice for the joint probability distribution for the random vector $(J_1^t,J_2^t,J_3^t)$	 implies that the random variables $J_1^t$ and $J_2^t$ are dependent and that $J_3^t$ is independent. It is easy to see that all these random variable have, separately, a Bernoulli distribution with parameter $1/2$, i.e., at every time-step one particle arrives at every node with probability $1/2$ and no particles arrive at every node with probability $1/2$.   In figure~\ref{fig:three-vertex-exact-sol} we show a graphical representation of the open chain. 

\begin{figure}[ht]
\begin{center}
\scalebox{0.5}{\includegraphics{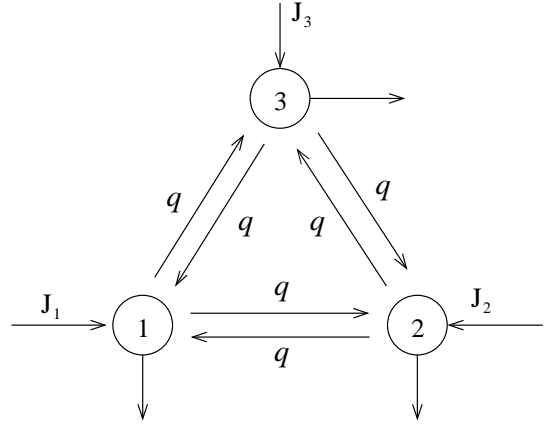}}
\end{center}
     \caption{
     An example of a three-states open Markov chain. 
     }
\label{fig:three-vertex-exact-sol}
\end{figure}
%

A straightforward calculation shows that the first and second moments of $\mathbf{J}^t$ are given by 
\begin{eqnarray}
\mathbb{E}[\mathbf{J}^t] &=& \bepsilon = (1/2, 1/2,1/2),
\\
\mathrm{Var}(\mathbf{J}^t) &=& \bDelta = \frac{1}{4} \left[
  \begin{array}{ccc}
  1 & 2p-1 & 0 \\ 
  2p-1 & 1 & 0 \\ 
  0 & 0 & 1
\end{array} \right].
\end{eqnarray}

Now we calculate the mean stationary distribution of particles, $\overline{\bmu}$, and the stationary variance matrix, $\overline{\bSigma}$. The calculation of $\overline{\bmu}$  is straightforward. We only need to obtain the inverse of the matrix $ \mathbf{1} -\mathbf{Q}$, which is given by
\[
( \mathbf{1} -\mathbf{Q} )^{-1} =  \frac{1}{ 1-3q^2-2q^3} \left[
  \begin{array}{ccc}
  1-q^2 & q+q^2 & q+q^2 \\ 
  q+q^2 & 1-q^2 & q+q^2 \\ 
  q+q^2 & q+q^2 & 1-q^2
\end{array} \right].
\]
With the above result we can see that the mean stationary distribution $\overline{\bmu}$ is given by
\[
\overline{\bmu} = \bepsilon ( \mathbf{1} -\mathbf{Q} )^{-1}  =  \left( \frac{1}{2-4q},\frac{1}{2-4q},\frac{1}{2-4q} \right).
\]
The last result implies that the  particles on the state space are equidistributed when the system has reached the stationarity. Moreover, we can also appreciate that the number of particles on every state diverges as the parameter $p$ tends to $1/2$. This divergence is actually a consequence of the fact that the system does not allow that the particles leave the state space when $p=1/2$. Thus, for such a parameter value,  the system is still receiving particles but it does not allow that the particles escape to the outside, thus increasing indefinitely the number of particles inside the system.

Now let us compute the stationary variance matrix $\overline{\bSigma}$. Recall that this quantity can be obtained by means of the formula,
\begin{equation}
\overline{\bSigma} = \sum_{k=0}^\infty (\mathbf{Q}^T)^k (\bDelta + \overline{\bLambda}) \mathbf{Q}^k.
\label{eq:bSigma-int}
\end{equation}
First let us obtain the explicit form of $\overline{\bLambda}$. According to equation~(\ref{eq:barLambda}) we have that 
\begin{eqnarray}
( \overline{\bLambda} )_{i,j} &=& \sum_{k=1}^S (\overline{\bmu})_k \left(q_{k,i}\delta_{i,j}-q_{k,i}q_{k,j}  \right)
\nonumber
\\
&=&  \frac{1}{2-4q} \sum_{k=1}^3 \left(q_{k,i}\delta_{i,j}-q_{k,i}q_{k,j}  \right)
\nonumber
\\
&=&  \frac{q}{1-2q}  \delta_{i,j} -  \frac{1}{2-4q} \sum_{k=1}^S q_{k,i}q_{k,j}.
\nonumber
\end{eqnarray}
Thus, we have
\[
\overline{\bLambda} = \frac{1}{2-4q}  \left[
  \begin{array}{ccc}
  2q(1-q) & -q^2 & -q^2 \\ 
  -q^2 & 2q(1-q) & -q^2 \\ 
  -q^2 & -q^2 & 2q(1-q)
\end{array} \right].
\]
Next we need to compute the $n$th power of the matrix $\mathbf{Q}$. It is not hard to prove that 
\begin{equation}
\label{eq:Q^n}
 \mathbf{Q}^k = \frac{q^n}{3} \left[
  \begin{array}{ccc}
  2^k + 2(-1)^k  & 2^n -(-1)^n & 2^k -(-1)^k \\ 
2^k -(-1)^n & 2^k + 2(-1)^n & 2^k -(-1)^k \\ 
2^k -(-1)^n & 2^k -(-1)^n & 2^k + 2(-1)^k
\end{array} \right].
\end{equation}
This result shows that for this particular case 
\[
(\mathbf{Q}^T)^k = \mathbf{Q}^k,
\]
because $\mathbf{Q}$ is itself a symmetric matrix. 

It is not hard to see that the term 
\[
\mathbf{Q}^k (\bLambda + \bDelta) \left(\mathbf{Q}^T\right)^k 
\]
is a matrix whose components are all exponential (or linear combinations of exponentials) in the variable $k$. This observation allows us to see that the infinite summation~(\ref{eq:bSigma-int}) can be exactly computed. Then, we can obtain a closed expression for $\overline{\bSigma}$ by using symbolic calculations performed in the software \emph{Mathematica}. Thus we obtain,
\begin{eqnarray}
\overline{\bSigma} &=&\left[
\begin{array}{ccc}
 \frac{-8 q^5+(8 p-2) q^4+4 q^3+3 q^2+4 q+1}{4 \left(8 q^6-6 q^4-3 q^2+1\right)} & \frac{2 q^4+q^2+p \left(-8 q^4-4
   q^2+2\right)-1}{4 \left(8 q^6-6 q^4-3 q^2+1\right)} & -\frac{q^2 \left(q^2-p+1\right)}{2 \left(8 q^6-6 q^4-3 q^2+1\right)} \\
 \frac{2 q^4+q^2+p \left(-8 q^4-4 q^2+2\right)-1}{4 \left(8 q^6-6 q^4-3 q^2+1\right)} & \frac{-8 q^5+(8 p-2) q^4+4 q^3+3 q^2+4
   q+1}{4 \left(8 q^6-6 q^4-3 q^2+1\right)} & -\frac{q^2 \left(q^2-p+1\right)}{2 \left(8 q^6-6 q^4-3 q^2+1\right)} \\
 -\frac{q^2 \left(q^2-p+1\right)}{2 \left(8 q^6-6 q^4-3 q^2+1\right)} & -\frac{q^2 \left(q^2-p+1\right)}{2 \left(8 q^6-6 q^4-3
   q^2+1\right)} & \frac{-8 q^5+(6-8 p) q^4+4 q^3+(4 p+1) q^2+4 q+1}{4 \left(8 q^6-6 q^4-3 q^2+1\right)} \\
\end{array}
\right]
\nonumber
\end{eqnarray}

Now, let us define the \emph{space correlation functions}. We will denote by $\kappa_{i,j}$ the correlation function between the random variables $N_i^t$ and $N_j^t$ as follows,
\begin{eqnarray}
\kappa_{i,j}  &:=& \mathrm{Corr}(N_i^t,N_j^t) =  \frac{\overline{\Sigma}_{i,j}}{\sqrt{\overline{\Sigma}_{i,i} \overline{\Sigma}_{j,j} }},
\end{eqnarray}
for all $i,j \in \mathcal{S}$. 
It is not hard to see that the correlation functions for this example are given by,
\begin{eqnarray}
\kappa_{1,2} &=&\frac{p \left(-8 q^4-4 q^2+2\right)+2 q^4+q^2-1}{\sqrt{\left((8 p-2) q^4-8 q^5+4 q^3+3 q^2+4 q+1\right)^2}},
\label{eq:k12}
\\
\kappa_{1,3} &=& -\frac{2 q^2 \left(-p+q^2+1\right)}{\sqrt{\left( (6-8 p) q^4+(4 p+1) q^2-8 q^5+4 q^3+4 q+1\right) }}
\nonumber
\\
&\times& \frac{1}{\sqrt{\left((8 p-2) q^4-8q^5+4 q^3+3 q^2+4 q+1\right)} },
\label{eq:k13}
\\
\kappa_{2,3} &=& \frac{-2 q^2 \left(-p+q^2+1\right)}{\sqrt{\left((6-8 p) q^4+(4 p+1) q^2-8 q^5+4 q^3+4 q+1\right) }}
\nonumber
\\
&\times& \frac{1}{\sqrt{\left((8 p-2) q^4-8q^5+4 q^3+3 q^2+4 q+1\right)} }.
\label{eq:k23}
\end{eqnarray}  

%
\begin{figure}[ht]
\begin{center}
\scalebox{0.6}{\includegraphics{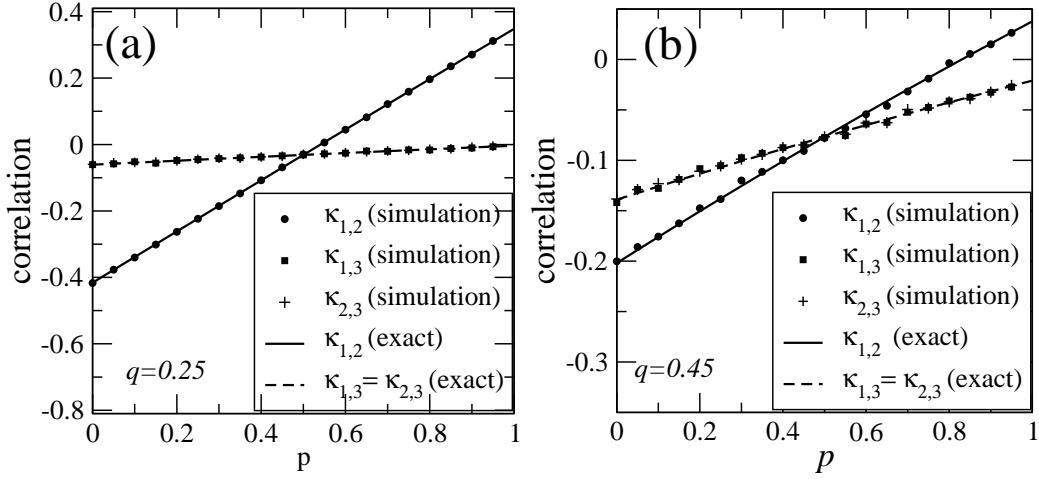}}
\end{center}
     \caption{
The correlation functions $\kappa_{1,2}$ and $\kappa_{2,3}$. In panel (a) we plot the correlation functions for the parameter value $q=0.25$ fixed and varying the parameter $p$. The solid line corresponds to $\kappa_{1,2}$ computed from equation~(\ref{eq:k12}) and the filled circles corresponds to  $\kappa_{1,2}$ numerically obtained from the simulations of the stochastic dynamics during $5\times 10^5$ time steps. The dashed line corresponds to $\kappa_{1,3} = \kappa_{2,3}$ computed from equation~(\ref{eq:k13}) and the filled squares corresponds to  $\kappa_{1,3}$ numerically obtained from the simulations 
. In panel (b) we do the same as in panel (a) but using the parameter value $q=0.45$. 
     }
\label{fig:spatial_correlations}
\end{figure}
%
In order to test the results we obtained for this example we performed numerical simulations. We have simulated the dynamics of the open chain for several parameter values during a time of $5\times 10^5$ steps. Then, from the time series obtained we estimated the correlation functions $\kappa_{1,2} $,  $\kappa_{1,3} $  and $\kappa_{2,3}$ and we compare them against the theoretical prediction given in equations~(\ref{eq:k12}), (\ref{eq:k13}) and~(\ref{eq:k23}) (note that for our example $\kappa_{1,3} = \kappa_{2,3}$). In figures~\ref{fig:spatial_correlations}a and ~\ref{fig:spatial_correlations}b we show the correlation functions obtained by numerical simulations and computed from  equations~(\ref{eq:k12}) and~(\ref{eq:k23}). We take the parameter value $q=0.25$ and plotted $\kappa_{1,2}$ (solid line, filled circles) and $\kappa_{2,3}$  (dashed line, filled squares) as a function of $p$ (figure~\ref{fig:spatial_correlations}a). The same graph is done but using the parameter value $q=0.45$ (figure~\ref{fig:spatial_correlations}b). As we can see, our the results obtained from numerical simulations are consistent with the formulas that we obtained theoretically.

\end{example}


\subsection{Distribution of particles leaving the state space}


Up to now we have given an expression for the m.g.f.~for the number of particles in the state space. Since the system is open, at every time step there is a number of particles arriving to the system, which is determined by the random vector $\mathbf{J}^t$. The total number of particles per step arriving to the system is then
\[
I_t = \sum_{i=1}^s J^t_i.
\]
When the system reaches the stationarity the mean number of particles within the system attains a constant value, meaning that there is a equilibration between the number of incoming and outcoming particles. To be precise, if we denote by $O_t$ the total number of particles leaving the state space, then, under stationarity we should have that
\[
\mathbb{E}[I_t] = \mathbb{E}[O_t],
\]
a fact that can be inferred by the ``conservation'' of the mean number of particles at stationarity. Our goal here is to go beyond the above expression, we would like to characterize how the random variable $O_t$ evolves in time and how much it is influenced by the incoming number of particles and the ``jumping'' rules of the Markov chain. To this end, let us define some quantities which will allow to describe the random variable $O_t$ explicitly.

\begin{definition}
\label{def:outcoming}

Let $\big( \mathcal{S}, \mathbf{Q}, \{\mathbf{J}^t : t \in \mathbb{N}\} \big)$ be an open Markov chain with state space $ \mathcal{S}$,  jump matrix $\mathbf{Q}$ (with components $(\mathbf{Q})_{i,j} = q_{i,j}$) and incoming protocol $\{ \mathbf{J}^{t}  : t \in \mathbb{N} \}$. Let $ e_i $ be defined as the escape probability from the state $i$, i.e.~the probability with which a particle in the state $i$ leaves the system to the outside,  
\begin{equation}
e_i := 1 -\sum_{j=1}^S q_{i,j}. 
\end{equation}
Next, let $\mathbf{U}^t = (U_1^t,U_2^t,\dots,U_S^t)$ be a random vector whose components have binomial distribution as follows,
\[
U_i^t \sim \mathrm{Binom}(e_i,N_i^t)
\]
where $\mathbf{N}^t$  is time-dependent distribution over the state space. Then we say the $U_i^t$ is the number of particles leaving the state $i$ to the outside at time $t$. The total number of particles $O_t$ leaving the system is then the random variable given by
\begin{equation}
O_t := \sum_{i=1}^S U_{i}^t.
\end{equation}

Finally, we define the vector $\mathbf{e}$, with components $e_i = (\mathbf{e})_i$, which will be referred to as the \emph{escape probability vector} of the chain. Additionally let $\mathbf{E} $ be  a diagonal  matrix with components $(\mathbf{E})_{i,j} = E_{i,j}$ defined as
\begin{equation}
E_{i,j}:= e_i \delta_{i,j},
\end{equation}
which will be referred to as the e \emph{escape probability matrix} of the chain.

\end{definition}

Our main goal is now to characterize the random vector $\mathbf{U}^t$ giving the number of particles leaving the chain. We should notice that the distribution of $ \mathbf{U}^t $ depends on $\mathbf{N}^t$ which is also a random vector. This implies that, given the value of $\mathbf{N}^t$, we can specify  the conditional distribution for $\mathbf{U}^t$, i.e., 
\begin{equation}
T(\mathbf{n};\mathbf{m}) := \mathbb{P}\left( \mathbf{U}^t = \mathbf{m} | \mathbf{N}^t = \mathbf{n}   \right),
\end{equation}
which is given by
\begin{equation}
T(\mathbf{n};\mathbf{m}) = \prod_{i=1}^S \frac{n_i!}{(n_i-m_i)!m_i!} e_i^{m_i} (1-e_i)^{n_i-m_i}.
\end{equation}
Once we know $T(\mathbf{n};\mathbf{m})$, we can write an expression for the probability distribution for the random vector $ \mathbf{U}^t$,
\begin{equation}
\label{eq:conditional1}
\mathbb{P} \left( \mathbf{U}^t=\mathbf{m} \right) = \sum_{ \mathbf{n}\in \mathbb{N}_0^S} \mathbb{P} \left(  \mathbf{U}^t = \mathbf{m} |  \mathbf{N}^t = \mathbf{n} \right) \mathbb{P} \left(    \mathbf{N}^t = \mathbf{n}  \right).
\end{equation}
If we denote by $ r_t (\mathbf{m})$ the probability distribution of the random vector $ \mathbf{U}^t$, i.e.,
\[
r_t(\mathbf{n}) := \mathbb{P} \left( \mathbf{U}^t=\mathbf{m} \right)
\]
then it is clear that equation~(\ref{eq:conditional1}) can be rewritten as
\begin{equation}
\label{eq:rt}
r_t(\mathbf{m}) = \sum_{n\in \mathbb{N}_0^S} p_t(\mathbf{n}) T(\mathbf{n};\mathbf{m}).
\end{equation}

The next step consists in obtaining the moment generating function of $\mathbf{U}^t$. This is because, as we saw above, the expression for the m.g.f.~of $\mathbf{N}_t$ can be written explicitly.  Thus, let $\mathcal{R}_t : \mathbb{R}^S \to \mathbb{R}$ be the m.g.f.~of $\mathbf{U}^t$, 
\begin{equation}
\mathcal{R}_t ( \balpha) := \mathbb{E} \left[ e^{\mathbf{U}^t \balpha^T} \right] = \sum_{\mathbf{m}\in \mathbb{N}_0^S} r_t(\mathbf{m}) e^{\mathbf{m}\balpha^T}. 
\end{equation}
Thus we can see that, using expression~(\ref{eq:rt}),  $\mathbf{R}_t$ can be written as follows 
\begin{equation}
\label{eq:Rt1}
\mathcal{R}_t (\balpha) = \sum_{\mathbf{m}\in \mathbb{N}_0^S}  \sum_{n\in \mathbb{N}_0^S} p_t(\mathbf{n}) T(\mathbf{n};\mathbf{m})  e^{\mathbf{m}\balpha^T}.
\end{equation}
Since $T (\mathbf{n};\mathbf{m})$ is a product of binomial distributions, it is clear that the sum over $\mathbf{m}$ results in the product of moment generating functions of binomial random variables,
\begin{equation}
\sum_{m\in\mathbb{N}_0^S}  T(\mathbf{n};\mathbf{m})  e^{\mathbf{m}\balpha^T} = \prod_{i=1}^S \left( 1-e_i + e_i e^{ \alpha_i}  \right)^{n_i}.
\end{equation}
Moreover, if we define the function $ \mathbf{C} = (C_1, C_2, \dots, C_S) : \mathbb{R}^S \to \mathbb{R}^S$ as follows,
\begin{equation}
C_i (\balpha ) = C_i (\alpha_i) = \log \left(  1-e_i + e_i e^{ \alpha_i} \right),  \quad \mathrm{for} \quad 1\leq i \leq S,
\end{equation}
it is clear that we can write 
\begin{equation}
\sum_{m\in\mathbb{N}_0^S}  T(\mathbf{n};\mathbf{m})  e^{\mathbf{m}\balpha^T} = \exp\left(\mathbf{n}\mathbf{C}^T(\balpha) \right).
\end{equation}
The above expression, together with equation~(\ref{eq:Rt1}), gives us
\begin{equation}
\mathcal{R}_t (\balpha ) = \mathcal{G}_t \left( \mathbf{C}(\balpha) \right).
\end{equation}
The last relation states that the moment generating function of $\mathbf{U}^t$ can be obtained by means of the moment generating function of $\mathbf{N}^t$, which is mediated by the transformation $\mathbf{C}$. 

Observe that $\mathcal{R}_t (\balpha )$ allows us to obtain the first and the second moment of number of leaving particles. Taking the first derivative (gradient) of $\mathcal{R}_t (\balpha )$ and evaluating it in $\balpha =\mathbf{0}$ we obtain,
\begin{equation}
\label{eq:meanU}
\mathbb{E}[\mathbf{U}^t] = \bmu_t \mathbf{E} = ( (\bmu_t)_1 e_1,  (\bmu_t)_2 e_2, \dots,  (\bmu_t)_S e_S).
\end{equation}
The above result allows us to calculate the mean number of leaving particles $\mathbb{E}[O_t]$,
\begin{equation}
\label{eq:mean_leaving}
\mathbb{E}[O_t] = \sum_{i=1}^S \mathbb{E}[U_i^t] =  \sum_{i=1}^S (\bmu_t)_i e_i = \bmu_t \mathbf{e}^T.
\end{equation}
Now, when the system reaches stationarity we have to replace $\bmu_t$ by $\overline{\bmu}$,
\[
(\mathbb{E}[O_t] )_{\mathrm{stat}} =  \sum_{i=1}^S (\overline{\bmu})_i e_i.
\]
Thus, from the definition of $e_i$ we obtain 
\begin{eqnarray}
(\mathbb{E}[O_t] )_{\mathrm{stat}} &=&  \sum_{i=1}^S (\overline{\bmu})_i \left( 1 - \sum_{j=1}^S  q_{i,j}\right)
\nonumber
\\
 &=&  \sum_{i=1}^S (\overline{ \bmu})_i  - \sum_{j=1}^S \sum_{i=1}^S  (\overline{\bmu})_i    q_{i,j} 
\nonumber
\\
 &=&  \sum_{i=1}^S (\overline{ \bmu})_i  - \sum_{j=1}^S  (\overline{\bmu} \mathbf{Q} )_j  = \sum_{i=1}^S \left( \overline{ \bmu} -  \overline{\bmu} \mathbf{Q} \right)_i,
\end{eqnarray}
and recalling that $\overline{\bmu}$ satisfy the equation
\[
\overline{\bmu} - \overline{\bmu} \mathbf{Q} = \bepsilon,
\]
we can observe that 
\begin{eqnarray}
(\mathbb{E}[O_t] )_{\mathrm{stat}} &=&   \sum_{i=1}^S (\bepsilon)_i =  \sum_{i=1}^S (\mathbb{E}[ \mathbf{J}^t])_i = (\mathbb{E}[I_t] )_{\mathrm{stat}},
\end{eqnarray}
which is the relation we have anticipated by invoking the particle number conservation principle.

The variance matrix $\mathrm{Var}(\mathbf{U}^t)$ can also be obtained by means of the second derivative of $\mathcal{R}_t (\balpha )$. A calculation achieved in~\ref{ape:1} shows that
\begin{equation}
\label{eq:varU}
\mathrm{Var}(\mathbf{U}^t) = \mathbf{E} \bSigma_t \mathbf{E} + \mathbf{D}_t,
\end{equation}
where $\mathbf{D}$ is a diagonal matrix with  components 
\[
(\mathbf{D}_t)_{i,j} :=  (\bmu_t)_i e_i (1-e_i)\delta_{i,j},
\]
where $\delta_{i,j}$ is the well-known Kr\"onecker delta. The variance matrix $\mathbf{U}^t$ allows us in turn to obtain a closed formula for the variance of the total number of leaving particles, 
\begin{eqnarray}
\label{eq:var_leaving}
\mathrm{Var}(O_t) &=& \sum_{i=1}^S  \sum_{j=1}^S  \left( \mathrm{Var}(\mathbf{U}^t) \right)_{i,j} 
\nonumber
\\
&=& \sum_{i=1}^S  \sum_{j=1}^S e_i ( \bSigma_t )_{i,j} e_j +   \sum_{i=1}^S  (\bmu_t)_i  (1-e_i) e_i
\\
&=& \mathbf{e} \bSigma_t \mathbf{e}^T + \bmu_t \left( \mathbf{1} - \mathbf{E}\right) \mathbf{e}^T,
\end{eqnarray}
where $\mathbf{1}$ stands for the $ S\times S$ identity matrix. We observe that the variance matrix of the number of outcoming particles is not the same as the variance matrix of the incoming particles nor the variance matrix of the particles in the system. Thus we have that the fluctuations and correlations are modulated by the internal dynamics of the system. This is important because measuring the correlations and fluctuations of the outcoming flux gives information on the internal dynamics of the system.

\begin{example}

Let us consider again the one-vertex model given in example~\ref{ex:one-vertex}. We have seen that the m.g.f.~of $\mathbf{N}^t$ is given by
\begin{eqnarray}
\mathcal{G}_{\mathrm{stat}} (\alpha)  &=&  \prod_{r=0}^{\infty}  \left( 1 - p q^r + p q^r e^{\alpha} \right).
\end{eqnarray}
Notice that the probability escape ``vector'' consists of a single number, given by $e_0 = 1-q$ (recall that the jump matrix is also a single number).  Thus we have that the transformation $C : \mathbb{R}  \to \mathbb{R}$ can be written as
\[
C(\alpha ) = \log \left( 1-e_0 + e_0 e^{\alpha} \right) = \log \left( q + (1-q) e^{\alpha} \right) .
\] 
Now we can obtain the m.g.f.~of $U$, the number of particles leaving the system,
\[
\mathcal{R}_{\mathrm{stat}} (\alpha ) = \mathcal{G}_{\mathrm{stat}} \left( C(\alpha)\right) = \prod_{r=0}^{\infty}  \left( 1 - p q^r + p q^r e^{C(\alpha)} \right).
\]
Notice that $e^{C(\alpha)}$ can be written as $e^{C(\alpha)} =  q + (1-q) e^{\alpha}$, thus we have, 
\begin{eqnarray}
\mathcal{R}_{\mathrm{stat}} (\alpha ) &=& \prod_{r=0}^{\infty}  \left( 1 - p q^r + p q^r \left( q + (1-q) e^{\alpha} \right) \right).
\nonumber
\\
&=& \prod_{r=0}^{\infty}  \left( 1 - p q^r (1-q) + p q^r (1-q) e^{\alpha}  \right).
\nonumber
\end{eqnarray}
The last expression establishes that the distribution of $U^t$ (at stationarity) can be seen as an infinite convolution of Bernoulli distributions, with parameters $ pq^r (1-q) $ for $r\in \mathbb{N}_0$. Thus, the random variable $U^t$  can be written as an infinite sum of i.i.d.~random variables $Y_r$ (with the above-mentioned Bernoulli distribution), 
\[
U^t = \sum_{r=0}^\infty Y_r.
\]
whenever $U^t$ has reached stationarity. Particularly the mean number of leaving particles as well as its variance can be exactly determined, 
\begin{eqnarray}
\mathbb{E}[U^t] &=&  \sum_{r=0}^\infty \mathbb{E}[Y^r] = \sum_{r=0}^\infty pq^r(1-q) = p.
\\
\mathrm{Var}(U^t) &=& \sum_{r=0}^\infty \mathrm{Var}(Y_r) = \sum_{r=0}^\infty pq^r(1-q) (1- pq^r(1-q))
\nonumber
\\
&=& p - \frac{p^2(1-q)^2}{1-q^2}.
\end{eqnarray}
We should emphasize that the above expressions can be obtained through the formulas for $\mathbb{E}[O_t]$ and  $\mathrm{Var}(O_t)$, given in equations~(\ref{eq:mean_leaving}) and~(\ref{eq:var_leaving}). Since $S=1$ (because we have only one state) we have that
\begin{eqnarray}
\mathbb{E}[O_t] &=& \mathbb{E}[N^t] e_0 = \frac{p}{1-q} e_0 = p,
\nonumber
\\
\mathrm{Var}(U^t) &=& \mathrm{Var}(N^t)e_0^2  +  \mathbb{E}[N^t] e_0 (1-e_0).
\nonumber
\\
&=& \left( \frac{p}{1-q} + \frac{p^2}{1-q^2}\right) (1-q)^2  + \left( \frac{p}{1-q} \right) q (1-q)
\nonumber
\\
&=& p + \frac{p^2(1-q)^2}{1-q^2}.
\end{eqnarray}
where we used the expressions for $\mathbb{E}[N^t]$ and  $\mathrm{Var}(N^t)$ given in equations~(\ref{eq:1v_mean}) and ~(\ref{eq:1v_var}).

\end{example}

\section{\label{sec:correlations}The influence of incoming particles on time-correlations}


\subsection{Time correlations for the open Markov chain}

Up to now we have seen that it is possible to find an explicit expression for the time-dependent distribution on the state space. This evolution is fully characterized by the two first cumulants, the mean distribution over the state space $\mathbf{\bmu}_t$ and the variance matrix $\mathbf{\bSigma}_t$. It is important to emphasize that we made no assumptions on the time correlations of the sequence of random vectors $\{ \mathbf{J}^t  : t\in \mathbb{N}\}$. This is because to obtain the distribution over the state space, $\mathbf{N}^t$, for a given time $t$, it was enough to know the number of incoming particles at time $t$. On the other hand, if we would like to compute the two-times correlation function for certain observable we necessarily have to known the number of incoming particles a two different times. This information unavoidably will be related to the two-times covariance matrix of the process $\{ \mathbf{J}^t  : t\in \mathbb{N}\}$, i.e., the covariance matrix between the random vector $\mathbf{J}^t$ and $\mathbf{J}^{t+s}$ for $s,t \in \mathbb{N}$. Our goal in this section is to obtain an expression for the covariance ${C}_{i,j}(t,t+s)$ between the $i$th coordinate of $\mathbf{N}^{t}$ and the $j$th coordinate of $\mathbf{N}^{t+s}$, i.e., 
\begin{equation}
\label{eq:corr-def}
{C}_{i,j}(t,t+s) := \mathbb{E}[N^{t}_i N^{t+s}_j] - \overline{\mu}_i \overline{\mu}_j,
\end{equation}
where $\overline{\mu}_i$ is the $i$th coordinate of the mean stationary distribution.

In order to compute the expected value $\mathbb{E}[N^{t}_i N^{t+s}_j]$ it is necessary to have an expression for the stationary joint distribution $P_{t,t+s} (\mathbf{n},\mathbf{m})$. This quantity is defined as,
\begin{equation}
P_{t,t+s} (\mathbf{n},\mathbf{m}) := \mathbb{P}\left(\mathbf{N}^{t}=\mathbf{n};\mathbf{N}^{t+s} = \mathbf{m} \right),
\end{equation}
Our goal here is to establish a method to obtain the joint distribution $P_{t} (\mathbf{n},\mathbf{m})$. Actually, we will first determine an expression for a more  general quantity. Let $P (\mathbf{n}_0,\mathbf{n}_1,\mathbf{n}_s)$ denote the joint probability function of the random vectors $\mathbf{N}^t, \mathbf{N}^{t+1},\dots,\mathbf{N}^{t+s}$~\footnote{Notice that, to be  strict, the probability function $P$ depends on $t,t+1,\dots, t+s$, and we should denote this dependence explicitly by using subscripts, i.e., $P = P_{t,t+1,\dots,t+1}$. However, we will not use such a notation by the sake of simplicity in further calculations. The same convention will be adopted for other ``multiple-times'' joint probability functions or its corresponding moment generating functions.},
\begin{equation}
\label{eq:multiNdef}
P (\mathbf{n}_0,\mathbf{n}_1,\dots,\mathbf{n}_s) := \mathbb{P} \left( \mathbf{N}^t = \mathbf{n}_0;\mathbf{N}^{t+1} = \mathbf{n}_1;\dots;\mathbf{N}^{t+s} = \mathbf{n}_s \right).
\end{equation}

First of all, let us introduce some notation that will be useful to perform  further calculations. Let $ f (\mathbf{j}_0,\mathbf{j}_1,\dots, \mathbf{j}_{s-1})$ be the joint probability function of the random vectors $\mathbf{J}^{t},\mathbf{J}^{t+1},\dots,\mathbf{J}^{t+s-1}$, i.e.,
\begin{equation}
 f(\mathbf{j}_0,\mathbf{j}_1,\dots, \mathbf{j}_{s-1}) := \mathbb{P} \left( \mathbf{J}^{t}=\mathbf{j}_0;\mathbf{J}^{t+1}=\mathbf{j}_1;\dots,\mathbf{J}^{t+s-1} =\mathbf{j}_{s-1}\right).
\end{equation}
Let us also denote by $h (\mathbf{r}; \mathbf{k})$ the probability function of the random vector $\mathbf{R}^t$, which, as we saw in section~\ref{sec:evolution}, depends on the value taken by the random vector $\mathbf{N}^{t}$ (a value which we denote by $\mathbf{k}$ in the probability function $h$). 

With the above-introduced notation it is possible to write the joint probability function, given in equation~(\ref{eq:multiNdef}), in terms of the probability functions $f$ and $h$, 
\begin{eqnarray}
P (\mathbf{n}_0,\mathbf{n}_1,\dots,\mathbf{n}_s)  &=&  \mathbb{P} \left( \mathbf{N}^t = \mathbf{n}_0;\mathbf{N}^{t+1} = \mathbf{n}_1;\dots;\mathbf{N}^{t+s} = \mathbf{n}_s \right)
\nonumber
\\
&=& \mathbb{P} \left( \mathbf{N}^t = \mathbf{n}_0;\mathbf{N}^{t+1} = \mathbf{n}_1;\dots;\mathbf{N}^{t+s} = \mathbf{n}_s \right)
\nonumber
\\
\nonumber
&=& \mathbb{P} \left(\mathbf{N}^{t+1} = \mathbf{n}_1;\mathbf{N}^{t+2} = \mathbf{n}_2;\dots;\mathbf{N}^{t+s} = \mathbf{n}_s \big|  \mathbf{N}^t = \mathbf{n}_0 \right) \mathbb{P}(\mathbf{N}^t=\mathbf{n}_0)
\\
\nonumber
&=& \mathbb{P} \left(\mathbf{J}^{t}+\mathbf{R}^{t} = \mathbf{n}_1;\dots;\mathbf{J}^{t+s-1}+\mathbf{R}^{t+s-1} = \mathbf{n}_s \big|  \mathbf{N}^t = \mathbf{n}_0 \right) p_{t}(\mathbf{n}_0).
\end{eqnarray}
Notice that the random vector $\mathbf{R}^{t+j}$ depends on $\mathbf{n}_{j}$ for  $1\leq j \leq s-1$, values which are given a priori. Thus,  the  random vectors $\mathbf{R}^{t},\mathbf{R}^{t+1},\dots,\mathbf{R}^{t+s-1}$ are all independent (because the values taken by the random vectors $\mathbf{N}^{t+j}$ for $0\leq j\leq s$ are all fixed), which allows us to write 
\begin{eqnarray}
\nonumber
& & P(\mathbf{n}_0,\mathbf{n}_1,\dots,\mathbf{n}_s)  =  \mathbb{P} \left( \mathbf{N}^t = \mathbf{n}_0;\mathbf{N}^{t+1} = \mathbf{n}_1;\dots;\mathbf{N}^{t+s} = \mathbf{n}_s \right)
\\
\nonumber
&=& \sum_{\mathbf{j}_0+\mathbf{r}_0=\mathbf{n}_1} \sum_{\mathbf{j}_1+\mathbf{r}_1=\mathbf{n}_2} \dots \sum_{\mathbf{j}_{s-1}+\mathbf{r}_{s-1}=\mathbf{n}_s} 
\mathbb{P} \left(\mathbf{J}^{t} = \mathbf{j}_0;\mathbf{J}^{t+1} = \mathbf{j}_1;\dots;\mathbf{J}^{t+s-1} = \mathbf{j}_{s-1} \right) 
\\
\nonumber
&\times& \mathbb{P} \left( \mathbf{R}^{t} = \mathbf{r}_0 \right)\mathbb{P} \left( \mathbf{R}^{t+1} = \mathbf{r}_1 \right) \dots \mathbb{P} \left( \mathbf{R}^{t+s-1} = \mathbf{r}_{s-1}\right)
p_{t}(\mathbf{n}_0).
\end{eqnarray}
In terms of the probability functions defined above we have,
\begin{eqnarray}
\nonumber
P(\mathbf{n}_0,\mathbf{n}_1,\dots,\mathbf{n}_s)  &=& 
\sum_{\mathbf{j}_0+\mathbf{r}_0=\mathbf{n}_1} \sum_{\mathbf{j}_1+\mathbf{r}_1=\mathbf{n}_2} \dots \sum_{\mathbf{j}_{s-1}+\mathbf{r}_{s-1}=\mathbf{n}_s} 
f (\mathbf{j}_0,\mathbf{j}_1,\dots, \mathbf{j}_{s-1})
\\
&\times& h (\mathbf{r}_0; \mathbf{n}_0) h (\mathbf{r}_1; \mathbf{n}_1) \dots h (\mathbf{r}_{s-1}; \mathbf{n}_{s-1})
p_{t}(\mathbf{n}_0).
\label{eq:multiP}
\end{eqnarray}
The above expression states that the joint distribution $P(\mathbf{n}_0,\mathbf{n}_1,\dots,\mathbf{n}_s) $ can be written in terms of the stationary distribution $p_{\mathrm{stat}}$ and the probability functions of the random vector $\mathbf{R}^{t}$ and the joint distribution of the random vectors $\mathbf{J}^{t},\mathbf{J}^{t+1},\dots,\mathbf{J}^{t+s-1}$, distributions that are given a priori. Once knowing the joint distribution $P(\mathbf{n}_0,\mathbf{n}_1,\dots,\mathbf{n}_s) $, we can compute the two-times joint distribution $P_{t,t+s}(\mathbf{n},\mathbf{m})$, 
\begin{equation}
P_{t,t+s}(\mathbf{n},\mathbf{m}) = \sum_{\mathbf{n}_1\in\mathbb{N}_0}\sum_{\mathbf{n}_2\in\mathbb{N}_0} \dots \sum_{\mathbf{n}_{s-1}\in\mathbb{N}_0}P(\mathbf{n},\mathbf{n}_1,\dots,\mathbf{n}_{s-1},\mathbf{m}).
\end{equation}
The above expression can be used to obtain the moment generating function of $(\mathbf{N}^t,\mathbf{N}^{t+s})$ and then the corresponding two-times covariance matrix $C_{i,j}(t,t+s) = \left( \mathrm{Cov}(\mathbf{N}^t,\mathbf{N}^{t+s}) \right)_{i,j}$. Those calculations are performed in~\ref{ape:1}, here we only write down the result,
\begin{equation}
\label{eq:covNtNt+s}
\mathrm{Cov} (\mathbf{N}^t, \mathbf{N}^{t+s}) =  \bSigma_t  \mathbf{Q}^{s}.
\end{equation}
Notice that equation~(\ref{eq:covNtNt+s}) is valid even if the system has not necessarily reached stationarity. If we assume that the system has attained the stationarity (which means that $\bSigma_t$ no longer depend on time), we obtain
\begin{equation}
\label{eq:covNtNt+s_ss}
\mathrm{Cov} (\mathbf{N}^t, \mathbf{N}^{t+s}) = \overline{ \bSigma}  \mathbf{Q}^{s}.
\end{equation}
Due to stationarity, it is clear that the covariance matrix depends only on $s$, the difference between the times $t$ and $t+s$. 
\begin{example}

Let us consider the system introduced in example~\ref{ex:3v_symmetric}. We should notice that we were able to obtain an exact expression for the stationary variance matrix $\overline{\bSigma}$. Thus, computing the covariance matrix involves only the product of two matrices. The resulting expression for the covariance matrix is too long to write down here. Thus, instead of giving explicitly the expression for the covariance matrix, we will show the theoretically computed time-dependent correlation functions defined as
\begin{equation}
\label{eq:time-correlation-function}
\tilde{C}_{i,j}(s) := \frac{ \left( \mathrm{Cov} (\mathbf{N}^t, \mathbf{N}^{t+s}) \right)_{i,j} }{\sqrt{ \overline{\Sigma}_{i,i} \overline{\Sigma}_{j,j}  }}
\end{equation}
%
\begin{figure}[ht]
\begin{center}
\scalebox{0.4}{\includegraphics{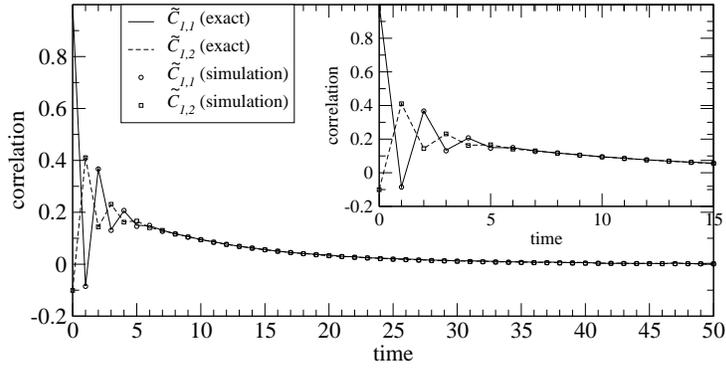}}
\end{center}
     \caption{
Time-dependent correlation function for the three-states open chain. We show the correlation functions for the parameter values $p=0.40$ and $q=0.45$. The theoretically computed correlation function  $\tilde{C}_{1,1}(s)$ is represented by the solid line and the numerically obtained from simulations are represented by the open cicles. Analogously, the theoretically computed correlation function  $\tilde{C}_{1,2}(s)$ is represented by the dashed line and the numerically obtained from simulations are represented by the open squares. The inset shows the same graph for small values of $s$.
     }
\label{fig:time-correlations}
\end{figure}
%

In figure~\ref{fig:time-correlations} we show the behavior of the correlation functions $\tilde{C}_{1,1}(s)$  and $\tilde{C}_{1,2}(s)$ for the three-states open Markov chain  studied in example~\ref{ex:3v_symmetric}. For the parameter values we display the theoretically computed  correlation function using equation~(\ref{eq:covNtNt+s_ss}) and~(\ref{eq:time-correlation-function}). The figure also shows the correlation functions obtained by means of numerical simulations of the system. The total time-steps performed to obtain the data from simulations was $5\times 10^5$. We appreciate that the theoretical results agree with the simulations showing the consistency of our results.

\end{example}

\section{Conclusions} 

We have introduced a simple model for open Markov chains by interpreting the state space of a usual Markov chain as physical ``sites'' where non-interacting particles can be placed and moving throughout it according to ``jumping rules'' given by a kind of stochastic matrix. The conditions for the chain to be open are given as a protocol of incoming particles, defined by a discrete-time stochastic process, and by a protocol of outcoming particles, which is implicitly defined by the condition that the ``stochastic matrix'' (called here jump matrix) has a spectral radius strictly less than one. These conditions establish the rules by means of which the particles arrive and leave the state space to the outside. We have shown that this model can be treated by means of the moment generating function technique, allowing us to obtain, in a closed form, the moment generating function of the distribution of particles over the state space. We have also shown that the system can be partially described by the dynamics of the two first cumulants of the the distribution of particles over the state space. Actually, we have given closed formulas for the two first cumulants when the system is able  to reach the stationarity. We have also studied how the correlations in the incoming protocol of particles are processed by the open chain. We have obtained closed formulas allowing compute the two-times covariance matrix for the random vector defined as the number of particles on the states.  Our main result is that the stationary two-times covariance matrix does not depend on the correlations of the particles arriving at the state space. This means that the stationary correlation functions essentially behaves as a closed Markov chain, i.e., that the correlations vanishes exponentially in time. The non-stationary correlations might probably content some information on the correlations of the incoming particles, but it would be necessary a more exhaustive study in this direction for a better understanding of such a process.

\begin{acknowledgments}
The author thanks the warm hospitality of the Instituto de F\'{\i}sica (IF-UASLP) and the Instituto Potosino de Investigaci\'on Cient\'{\i}fica y Tecnol\'ogica (IPICyT) where part of this work was achieved.
\end{acknowledgments}

\appendix

\section{\label{ape:1} Some mathematical calculations}

In this appendix we perform some calculations. Let us start by proving that the function $\mathcal{H} (\balpha)$ defined as
\begin{eqnarray}
\mathcal{H} (\balpha) &:=& \mathbb{E} \left[  e^{\mathbf{R}^t}\ba^T\right] 
= \sum_{\mathbf{r}\in\mathbb{N}_0^S }  \mathbb{P} (\mathbf{R}^t = \mathbf{r})  e^{\mathbf{r}\balpha^T }.
\end{eqnarray}
can be written as
\begin{equation}
\mathcal{H}(\ba) = e^{\mathbf{k} \mathbf{H}^T(\ba) },
\end{equation}
as it was stated in equation~(\ref{eq:gen-H-k}). First recall that the $j$th component of the random vector $\mathbf{R}^t$ is defined as, 
\begin{equation}
R^t_j = \sum_{i=1}^S B_{i,j}^t,
\end{equation}
where $\{ B_{i,j}^t : t \in \mathbb{N}_0, 1\leq i \leq S, 1\leq j \leq S \}$ is a sequence of i.i.d.~random variables with binomial distribution, $B_{ij}^t \sim \mathrm{Binom}( q_{i,j}, N^t_i)$ for all $t \in \mathbb{N}_0$, where $q_{i,j}$ is the $(i,j)$th component of $\mathbf{Q}$.  By independence it is clear that
\begin{eqnarray}
\mathbb{E} \left[  e^{\mathbf{R}^t\ba^T}\right]  = \prod_{j=1}^S \mathbb{E} \left[  e^{\alpha_j R_j^t}\right]  = \prod_{j=1}^S \prod_{i=1}^S \mathbb{E} \left[  e^{\alpha_j B_{i,j}^t}\right], 
\end{eqnarray}
Observe that the expected value 
\[
\mathbb{E} \left[  e^{\alpha_j B_{i,j}^t}\right]
\]
is the m.g.f.~of $B_{i,j}^t$ evaluated in $\alpha_j$. Thus we have
\begin{equation}
\label{eq:H1}
\mathcal{H} (\balpha) = \prod_{j=1}^S \prod_{i=1}^S \left( 1-q_{i,j} + q_{i,j}e^{\alpha_j} \right)^{k_i}.
\end{equation}
The above expression ca be rewritten, by convenience, as 
\begin{equation}
\mathcal{H} (\balpha) =  \prod_{i=1}^S \exp\left( k_i H_i(\balpha) \right) = \exp \left( \sum_{i=1}^S k_i H_i \right).
\end{equation}
where $H_i (\balpha)$ is defined as 
\[
H_i(\balpha):=\log \left(   \prod_{j=1}^S \left[ 1-q_{i,j} + q_{i,j}e^{\alpha_j} \right) \right) = \sum_{j=1}^S \log \left( 1-q_{i,j} + q_{i,j}e^{\alpha_j} \right).
\]
Thus, from equation~(\ref{eq:H1}), it is clear that the m.g.f.~$\mathcal{H}(\balpha)$ can be written as 
\begin{equation}
\mathcal{H}(\ba) = e^{\mathbf{k} \mathbf{H}^T(\ba) },
\end{equation}
where the (row) vector function $\mathbf{H}(\ba)$ is defined as
\[
\mathbf{H} (\balpha) = \left(H_1(\balpha),H_2(\balpha), \dots,H_S(\balpha)  \right).
\]

Next, we shall prove that the successive iterations of $ \mathbf{H}(\balpha)$ tend to $\mathbf{0} \in \mathbb{R}^S$. This can be shown by using the Contracting Principle in discrete dynamical systems and a lema on nonlinear maps in higher dimensions~\cite{katok1997introduction}. The only we need to do is to prove that $\mathbf{H}(\balpha)$ defines a contracting map in an open neighborhood around $\balpha = \mathbf{0}\in \mathbb{R}$, a point where the map $\mathbf{H}$ has a fixed point.

First notice that the  $\balpha=\mathbf{0}$ is a fixed point for the map $\mathbf{H} : \mathbb{R} \to \mathbb{R}$. Observe also that $\mathbf{H}$ is a continuously differentiable nonlinear map  whose first derivative (gradient) is given by
\begin{eqnarray}
\label{eq:1stderivative}
(D\mathbf{H})_{i,j} = \frac{\partial H_i }{\partial \alpha_j} = \frac{q_{i,j}e^{\alpha_j}}{ 1-q_{i,j} + q_{i,j}e^{\alpha_j}}.
\end{eqnarray}
If we evaluate $D\mathbf{H} (\balpha)$ at $\balpha=\mathbf{0}$ we have that 
\[
 \frac{ \partial H_i}{\partial \alpha_j }\bigg|_{\balpha = \mathbf{0}} = (\mathbf{Q})_{i,j} = q_{i,j},
\]
or, equivalently 
\[
D\mathbf{H} (\mathbf{0}) = \mathbf{Q}.
\]
Since the spectral radius of $\mathbf{Q} = D\mathbf{H} (\mathbf{0})$ is strictly less than one, then according to Lema 3.3.6 in Ref.~\cite{katok1997introduction}, we have that there is a closed neighborhood $U$ of $\balpha = 0$ such that $\mathbf{H}(U)\subset U$ where $\mathbf{H}$ is a eventually contracting map. Invoking the contracting principle for eventually contracting maps (cf. for example Corollary 2.6.13 in Ref.~\cite{katok1997introduction}) we have that, under iterates of $\mathbf{H}$, all points $\balpha \in U$ converge to $\mathbf{0}$. This proves that there is a neighborhood $U$ around $\balpha = \mathbf{0}$ such that $\mathbf{H}^{t}(\balpha) \to \mathbf{0}$ as $t\to \infty$ for all $\balpha \in U$.

Another relationship that it is important to prove here is
\begin{equation}
\frac{\partial^2 H_k (\balpha)}{ \partial \alpha_i \partial \alpha_j } \bigg|_{\balpha = \mathbf{0}} = q_{k,i} (1-q_{k,i})  \delta_{i,j}.
\end{equation}
which was used to obtain the dynamics of the second cumulant of $\mathbf{N}^t$ given in equation~(\ref{eq:varNt+1}). We have already computed the first derivative of $H_i$ with respect to $\alpha_j$ in equation~(\ref{eq:1stderivative}). Using this information we can compute the second derivative, obtaining 
\begin{eqnarray}
\frac{\partial^2 H_k (\balpha)}{ \partial \alpha_i \partial \alpha_j } \bigg|_{\balpha = \mathbf{0}} &=& 
\frac{\partial}{\partial \alpha_i} \left( \frac{q_{k,j} e^{\alpha_j}}{ 1-q_{k,j} + q_{k,j}e^{\alpha_j}} \right)
\bigg|_{\balpha = \mathbf{0}}
\nonumber
\\
& = &   \frac{  \left(  1-q_{k,j} + q_{k,j}e^{\alpha_j}  \right) q_{k,j} e^{\alpha_j} - \left( q_{k,j} e^{\alpha_j}  \right)^2   }{  \left(  1-q_{k,j} + q_{k,j}e^{\alpha_j} \right)^2 } \delta_{i,j} 
\Bigg|_{\balpha = \mathbf{0}}
\nonumber
\\
& = &q_{k,i} (1-q_{k,i})  \delta_{i,j},
\end{eqnarray}
which is the relationship we anticipated in equation~(\ref{eq:varNt+1}).

Next, we shall prove that the stationary solution to the cumulant dynamics equations are given by
\begin{eqnarray}
\overline{\bmu} &=& \bepsilon (\mathbf{1} - \mathbf{Q})^{-1}
\\
\overline{\bSigma} &=& \sum_{k=0}^\infty (\mathbf{Q}^T)^k  (\bDelta+\overline{\bLambda}  )\mathbf{Q}^k,
\end{eqnarray}
as it was stated in equations~(\ref{eq:mu_stat}) and~(\ref{eq:sigma_stat}) 

Let us start by recalling the equations for the dynamics of two first cumulants, given by expressions~(\ref{eq:mut+1}) and~(\ref{eq:varNt+1}), 
\begin{eqnarray}
\bmu_{t+1} &=& \bepsilon_t   + \bmu_t \mathbf{Q},
\label{eq:ape:mu_t+1}
\\
\bSigma_{t+1} &=&  \bDelta_t+\bLambda_t + \mathbf{Q}^T \bSigma_t \mathbf{Q},
\label{eq:ape:var_t+1}
\end{eqnarray}
where $\bLambda_t$ is defined in equation~(\ref{eq:Lambda_t}) as,
\begin{equation}
(\bLambda_t)_{i,j} := \sum_{k=1}^S (\bmu_t)_k q_{k,i}(1-q_{k,i}) \delta_{i,j}.
\end{equation}
The above means that $\bLambda_t$ is actually a function of $\bmu_t$, i.e., $\bLambda_t = \bLambda(\bmu_t)$, implying that the equations for the cumulant dynamic are not independent. Before going on the proof of the existence of a the stationary cumulants let us first impose the assumption, on equations~(\ref{eq:ape:mu_t+1}) and~(\ref{eq:ape:var_t+1}), that the incoming protocols is a stationary process. Assuming stationarity of the process $\{\mathbf{J}^t : t\in \mathbb{N}\}$ implies that the mean value of $\mathbf{J}^t$, as well as its variance, are independent on time. Thus, we can substitute $\bepsilon_t$ and $\bDelta_t$ by  $\bepsilon$ and $\bDelta$, respectively. Therefore, the dynamic equations for the cumulants can be rewritten as
\begin{eqnarray}
\bmu_{t+1} &=& \bepsilon   + \bmu_t \mathbf{Q},
\label{eq:ape:mu_t+1_ss}
\\
\bSigma_{t+1} &=&  \bDelta+\bLambda (\bmu_t) + \mathbf{Q}^T \bSigma_t \mathbf{Q},
\label{eq:ape:var_t+1_ss}
\end{eqnarray}
It is easy to see that the above recurrence equations can be solved to obtain,
\begin{eqnarray}
\bmu_{t+1} &=& \bmu_0 \mathbf{Q}^t +  \sum_{r=0}^{t-1}\bepsilon \mathbf{Q}^r ,
\label{eq:ape:mu_t+1_sol}
\\
\bSigma_{t+1} &=& \left( \mathbf{Q}^T\right)^t \bSigma_0 \mathbf{Q}^t +  \sum_{r=0}^{t-1} \left( \mathbf{Q}^T\right)^r \left( \bDelta + \bLambda (\bmu_{t-r}) \right) \mathbf{Q}^r.
\label{eq:ape:var_t+1_sol}
\end{eqnarray}
Finally, we should observe that the matrix $\mathbf{Q}^t$ (and its transpose) vanishes as $t\to\infty$. This is because the spectral radius of $\mathbf{Q}$ is strictly less than one, implying that successive powers of $\mathbf{Q}$ converge to the operator zero. This fact implies that the terms $ \bmu_0 \mathbf{Q}^t$ and $\left( \mathbf{Q}^T\right)^t \bSigma_0 \mathbf{Q}^t$  tend to zero as $t$ approaches infinity. Then, taking the limit of $t\to\infty$, we can see that equations~(\ref{eq:ape:mu_t+1_sol}) and~(\ref{eq:ape:var_t+1_sol}) result in,
\begin{eqnarray}
\overline{\bmu}:=\lim_{t\to\infty}\bmu_{t+1} &=& \sum_{r=0}^{\infty}\bepsilon \mathbf{Q}^r,
\label{eq:ape:mu_ss}
\\
\overline{\bSigma}:= \lim_{t\to\infty}\bSigma_{t+1} &=&  \sum_{r=0}^{\infty} \left( \mathbf{Q}^T\right)^r \left( \bDelta + \bLambda (\overline{\bmu}) \right) \mathbf{Q}^r.
\label{eq:ape:var_ss}
\end{eqnarray}
A straightforward calculations shows that $\overline{\bmu}$ and $\overline{\bSigma}$ are invariant under the dynamics given by equations~(\ref{eq:ape:mu_t+1_ss}) and~(\ref{eq:ape:var_t+1_ss}). We can see that the above expressions for $\overline{\bmu}$ and $\overline{\bSigma}$ are the ones anticipated in equations~(\ref{eq:mu_stat}) and~(\ref{eq:sigma_stat}).

Now let us compute the variance matrix of the random vector $\mathbf{U}^t$. First we should recall that the m.g.f.~of $\mathbf{U}^t$ is given by 
\begin{equation}
\label{eq:ape:Rt}
\mathcal{R}_t (\balpha ) = \mathcal{G}_t \left( \mathbf{C}(\balpha) \right),
\end{equation}
where $\mathbf{C}: \mathbb{R}^S \to \mathbb{R}^S $ is a transformation defined as
\begin{equation}
\left( \mathbf{C}(\balpha ) \right)_i = C_i (\alpha_i) := \log \left(  1-e_i + e_i e^{ \alpha_i} \right),  \quad \mathrm{for} \quad 1\leq i \leq S. 
\end{equation}
It is clear that the m.g.f.~allows us to obtain the variance matrix as  follows,
\begin{eqnarray}
\left(\mathrm{Var}(\mathbf{U}^t) \right)_{i,j} = \frac{\partial^2 \mathcal{R}_t }{ \partial \alpha_i \partial \alpha_j }\bigg|_{\balpha =\mathbf{0}} - \mathbb{E} [ U_i] \mathbb{E} [ U_j] .
\end{eqnarray}
Using expression~(\ref{eq:ape:Rt}) and the fact that $\mathbb{E}[U_i]=(\bmu_t \mathbf{E})_i $ (see equation~(\ref{eq:meanU}))  we obtain for the variance matrix
\begin{eqnarray}
\left(\mathrm{Var}(\mathbf{U}^t) \right)_{i,j} &=& \frac{\partial^2 \mathcal{G}_t \left( \mathbf{C}(\balpha) \right) }{ \partial \alpha_i \partial \alpha_j }
\bigg|_{\balpha =\mathbf{0}} - \left( \bmu_t \mathbf{E} \right)_i  \left( \bmu_t \mathbf{E} \right)_j
\nonumber
\\
&=& \frac{\partial }{\partial \alpha_i } \left(  \frac{\partial \mathcal{G}}{\partial C_j } \frac{\partial C_j }{\partial \alpha_j  } \right)
\bigg|_{\balpha =\mathbf{0}} - \left( \bmu_t \mathbf{E} \right)_i  \left( \bmu_t \mathbf{E} \right)_j
\nonumber
\\
&=&
\left(  \frac{\partial^2 \mathcal{G}}{\partial C_i \partial C_j } \frac{\partial C_i }{\partial \alpha_i  }\frac{\partial C_j }{\partial \alpha_j  } 
+  \frac{\partial \mathcal{G}}{\partial C_j } \frac{\partial^2 C_j }{\partial \alpha_j^2  } \delta_{i,j}\right)
\bigg|_{\balpha =\mathbf{0}} - \left( \bmu_t \mathbf{E} \right)_i  \left( \bmu_t \mathbf{E} \right)_j.
\nonumber
\end{eqnarray}
At this point it is important recalling that 
\begin{eqnarray}
\frac{\partial \mathcal{G}}{\partial \alpha_i }  \bigg|_{\balpha =\mathbf{0}} &=& (\bmu_t)_i,
\nonumber
\\
\frac{\partial^2 \mathcal{G}}{\partial \alpha_i \partial \alpha_j }  \bigg|_{\balpha =\mathbf{0}} &=&  \mathbb{E}[N^t_i N^t_j],
\nonumber
\end{eqnarray}
Additionally, some elementary calculations allow us to observe that the derivatives of $\mathbf{C}$ result in
\begin{eqnarray}
\frac{\partial C_i}{\partial \alpha_i }  \bigg|_{\balpha =\mathbf{0}} &=& e_i,
\nonumber
\\
\frac{\partial^2 C_i}{\partial \alpha_i^2 }  \bigg|_{\balpha =\mathbf{0}} &=& e_i (1-e_i). 
\end{eqnarray}
The above results, together with the fact that $\mathbf{C} (\mathbf{0})  =\mathbf{0}$, imply that the variance matrix can be written as
\begin{eqnarray}
\left(\mathrm{Var}(\mathbf{U}^t) \right)_{i,j} &=& 
 \mathbb{E}[N^t_i N^t_j] e_i e_j 
+   (\bmu_t)_j e_j (1-e_j) \delta_{i,j} - \left( \bmu_t \mathbf{E} \right)_i  \left( \bmu_t \mathbf{E} \right)_j.
\nonumber
\\
&=&
e_i (\bSigma_t)_{i,j} e_j + (\bmu_t)_i (\bmu_t)_i e_ie_j +   (\bmu_t)_j e_j (1-e_j) \delta_{i,j} - \left( \bmu_t \mathbf{E} \right)_i  \left( \bmu_t \mathbf{E} \right)_j,
\nonumber
\end{eqnarray}
where we used the fact that  $\mathbb{E}[N^t_i N^t_j] e_i e_j = (\bSigma_t)_{i,j} e_j + (\bmu_t)_i (\bmu_t)_i $. Now recalling that that $(\mathbf{E})_{i,j} =  e_i \delta_{i,j}$ it is clear that  
\begin{eqnarray}
\left(\mathrm{Var}(\mathbf{U}^t) \right)_{i,j}  &=&
\left( \mathbf{E} \bSigma_t \mathbf{E} \right)_{i,j}  +   (\bmu_t)_j e_j (1-e_j) \delta_{i,j},
\end{eqnarray}
which is the expression anticipated in equation~(\ref{eq:varU}).

Next, our goal consists in obtaining the two-times covariance matrix $\mathrm{Cov}(\mathbf{N}^t,\mathbf{N}^{t+s})$. To this end we need to obtain an expression for the m.g.f.~of the joint distribution $P_{t,t+s}(\mathbf{n},\mathbf{m})$. We have seen that the latter can be expressed in terms of the multiple-times joint distribution as follows
\begin{equation}
P_{t,t+s}(\mathbf{n},\mathbf{m}) = \sum_{\mathbf{n}_1\in\mathbb{N}_0}\sum_{\mathbf{n}_2\in\mathbb{N}_0} \dots \sum_{\mathbf{n}_{s-1}\in\mathbb{N}_0}
P(\mathbf{n},\mathbf{n}_1,\dots,\mathbf{n}_{s-1},\mathbf{m}).
\end{equation}
Thus, the moment generating function $\mathcal{L}_{t,t+s}$ of the random vectors $(\mathbf{N}^t,\mathbf{N}^{t+s})$ can be written as

\begin{eqnarray}
\mathcal{L}_{t,t+s} (\balpha,\bbeta) &:=& \sum_{\mathbf{n}_0\in\mathbb{N}_0 }   \sum_{\mathbf{n}_s\in\mathbb{N}_0 } P_{t,t+s}(\mathbf{n}_0,\mathbf{n}_s) e^{\mathbf{n}_0 \balpha^T} e^{\mathbf{n}_s \bbeta^T}
\nonumber
\\
&=&
\sum_{\mathbf{n}_0\in\mathbb{N}_0 } \sum_{\mathbf{n}_1\in\mathbb{N}_0 } \cdots  \sum_{\mathbf{n}_s\in\mathbb{N}_0 } 
 P(\mathbf{n}_0,\mathbf{n}_1,\dots,\mathbf{n}_{s-1},\mathbf{n}_s)  e^{\mathbf{n}_0 \balpha^T} e^{\mathbf{n}_s \bbeta^T}. \quad
 \nonumber
\end{eqnarray}

Using expression~(\ref{eq:multiP}) for the joint distribution $ P(\mathbf{n}_0,\mathbf{n}_1,\dots,\mathbf{n}_{s-1},\mathbf{n}_s)$ we obtain,
\begin{eqnarray}
\mathcal{L}_{t,t+s} (\balpha,\bbeta) &:=& \sum_{\mathbf{n}_0\in\mathbb{N}_0 } \sum_{\mathbf{n}_1\in\mathbb{N}_0 } \cdots  \sum_{\mathbf{n}_s\in\mathbb{N}_0 } 
 \sum_{\mathbf{j}_0+\mathbf{r}_0=\mathbf{n}_1} \sum_{\mathbf{j}_1+\mathbf{r}_1=\mathbf{n}_2} \cdots \sum_{\mathbf{j}_{s-1}+\mathbf{r}_{s-1}=\mathbf{n}_s} 
 f (\mathbf{j}_0,\mathbf{j}_1,\dots, \mathbf{j}_{s-1})
 \nonumber
 \\
&\times&
 h (\mathbf{r}_0; \mathbf{n}_0) h (\mathbf{r}_1; \mathbf{n}_1) \dots h (\mathbf{r}_{s-1}; \mathbf{n}_{s-1})
p_{t}(\mathbf{n}_0)
 e^{\mathbf{n}_0 \balpha^T} e^{\mathbf{n}_s \bbeta^T}. \quad
 \nonumber
\end{eqnarray}
Observe that, as in previous calculations, the summation over $\mathbf{n}_{j+1}$ together with the restricted summation over $\mathbf{j}_j + \mathbf{r}_j = \mathbf{n}_{j+1}$ results in a double summations over two independent (unrestricted) indices $\mathbf{j}_j $ and $\mathbf{r}_j $ for $0\leq j < s$. This fact leads us to
\begin{eqnarray}
\mathcal{L}_{t,t+s} (\balpha,\bbeta) &:=& \sum_{\mathbf{n}_0\in\mathbb{N}_0 } 
\sum_{\mathbf{j}_0\in\mathbb{N}_0 }\sum_{\mathbf{j}_1\in\mathbb{N}_0 } \cdots  \sum_{\mathbf{j}_{s-1}\in\mathbb{N}_0 } 
\sum_{\mathbf{r}_0\in\mathbb{N}_0 }\sum_{\mathbf{r}_1\in\mathbb{N}_0 } \cdots  \sum_{\mathbf{r}_{s-1}\in\mathbb{N}_0 } 
 p_{t}(\mathbf{n}_0)  f (\mathbf{j}_0,\mathbf{j}_1,\dots, \mathbf{j}_{s-1})
 \nonumber
 \\
&\times&
 h (\mathbf{r}_0; \mathbf{n}_0) h (\mathbf{r}_1; \mathbf{j}_0 + \mathbf{r}_0 )  \dots h (\mathbf{r}_{s-1}; \mathbf{j}_{s-2} + \mathbf{r}_{s-2})
 e^{\mathbf{n}_0 \balpha^T} e^{\mathbf{j}_{s-1} \bbeta^T} e^{\mathbf{r}_{s-1} \bbeta^T}. \quad
 \nonumber
 \\
\label{eq:ape:inter-multi-mgf}
\end{eqnarray}
We should notice that in the above expression, the summation over $\mathbf{r}_{s-1}$ can be achieved because from the summand we can factorize a term that depends on $\mathbf{r}_{s-1}$. Such a term is indeed $ h (\mathbf{r}_{s-1}; \mathbf{j}_{s-2} + \mathbf{r}_{s-2}) e^{\mathbf{r}_{s-1} \bbeta^T}$. We should recall that  $h( \mathbf{r}; \mathbf{n})$ is the probability function of a sum of random vectors with binomial distribution. We have shown above that the corresponding m.g.f.~is given by
\[
\mathcal{H}(\bbeta) := \sum_{\mathbf{r}\in\mathbb{N}_0} h( \mathbf{r}; \mathbf{n}) e^{\mathbf{r} \bbeta^T} = e^{\mathbf{n}\mathbf{H}^T(\bbeta)}.
\] 
The above implies that the summation over $\mathbf{r}_{s-1}$ in equation~(\ref{eq:ape:inter-multi-mgf}) results in 
\[
\sum_{\mathbf{r}_{s-1}\in\mathbb{N}_0}  h (\mathbf{r}_{s-1}; \mathbf{j}_{s-2} + \mathbf{r}_{s-2}) e^{\mathbf{r}_{s-1} \bbeta^T} = 
e^{(\mathbf{j}_{s-2} + \mathbf{r}_{s-2} )\mathbf{H}^T(\bbeta)}.
\] 
We then obtain,
\begin{eqnarray}
\mathcal{L}_{t,t+s} (\balpha,\bbeta) &:=& \sum_{\mathbf{n}_0\in\mathbb{N}_0 } 
\sum_{\mathbf{j}_0\in\mathbb{N}_0 }\sum_{\mathbf{j}_1\in\mathbb{N}_0 } \cdots  \sum_{\mathbf{j}_{s-1}\in\mathbb{N}_0 } 
\sum_{\mathbf{r}_0\in\mathbb{N}_0 }\sum_{\mathbf{r}_1\in\mathbb{N}_0 } \cdots  \sum_{\mathbf{r}_{s-1}\in\mathbb{N}_0 } 
 p_{t}(\mathbf{n}_0)  f (\mathbf{j}_0,\mathbf{j}_1,\dots, \mathbf{j}_{s-1})
 \nonumber
 \\
&\times&
 h (\mathbf{r}_0; \mathbf{n}_0) h (\mathbf{r}_1; \mathbf{j}_0 + \mathbf{r}_0 )  \dots h (\mathbf{r}_{s-2}; \mathbf{j}_{s-3} + \mathbf{r}_{s-3})
 \exp\left( {\mathbf{n}_0 \balpha^T} \right)
 \nonumber
 \\
 &\times&
 \exp\left( {\mathbf{j}_{s-1} \bbeta^T}  +  \mathbf{j}_{s-2}\mathbf{H}^T(\bbeta) \right) \exp\left( \mathbf{r}_{s-2} \mathbf{H}^T(\bbeta) \right). \quad
\label{eq:ape:inter2-multi-mgf}
\end{eqnarray}
Next, in the above expression we can perform the summation over $\mathbf{r}_{s-2}$ in a similar way as we did it for $\mathbf{r}_{s-1}$. Actually, all the terms containing $\mathbf{r}_{s-2}$  are  $ h (\mathbf{r}_{s-2}; \mathbf{j}_{s-3} + \mathbf{r}_{s-3})$ and $ e^{\mathbf{r}_{s-2}  \mathbf{H}^T(\bbeta) }$. Using the same reasoning used above we see that
\[
\sum_{\mathbf{r}_{s-2}\in\mathbb{N}_0}  h (\mathbf{r}_{s-2}; \mathbf{j}_{s-3} + \mathbf{r}_{s-3}) e^{\mathbf{r}_{s-2}  \mathbf{H}^T(\bbeta)} = 
\exp\left( {(\mathbf{j}_{s-3} + \mathbf{r}_{s-2} )\left( \mathbf{H}^{(2)}(\bbeta)\right)^{T}}\right).
\] 
We then obtain,
\begin{eqnarray}
\mathcal{L}_{t,t+s} (\balpha,\bbeta) &=& \sum_{\mathbf{n}_0\in\mathbb{N}_0 } 
\sum_{\mathbf{j}_0\in\mathbb{N}_0 }\sum_{\mathbf{j}_1\in\mathbb{N}_0 } \cdots  \sum_{\mathbf{j}_{s-1}\in\mathbb{N}_0 } 
\sum_{\mathbf{r}_0\in\mathbb{N}_0 }\sum_{\mathbf{r}_1\in\mathbb{N}_0 } \cdots  \sum_{\mathbf{r}_{s-1}\in\mathbb{N}_0 } 
 p_{t}(\mathbf{n}_0)  f (\mathbf{j}_0,\mathbf{j}_1,\dots, \mathbf{j}_{s-1})
 \nonumber
 \\
&\times&
 h (\mathbf{r}_0; \mathbf{n}_0) h (\mathbf{r}_1; \mathbf{j}_0 + \mathbf{r}_0 )  \dots h (\mathbf{r}_{s-3}; \mathbf{j}_{s-4} + \mathbf{r}_{s-4})
 \exp\left( {\mathbf{n}_0 \balpha^T} \right)
 \nonumber
 \\
 &\times&
 \exp\left( {\mathbf{j}_{s-1} \bbeta^T}  +  \mathbf{j}_{s-2}\mathbf{H}^T(\bbeta)  + \mathbf{j}_{s-3} \left( \mathbf{H}^{(2)}(\bbeta)\right)^{T}\right)
 \exp\left( \mathbf{r}_{s-3} \left( \mathbf{H}^{(2)}(\bbeta)\right)^{T} \right). \quad
 \nonumber
 \\
\label{eq:ape:inter3-multi-mgf}
\end{eqnarray}
Proceeding inductively it is clear that the above expression results in
\begin{eqnarray}
\mathcal{L}_{t,t+s} (\balpha,\bbeta) &=& \sum_{\mathbf{n}_0\in\mathbb{N}_0 } 
\sum_{\mathbf{j}_0\in\mathbb{N}_0 } \sum_{\mathbf{j}_1\in\mathbb{N}_0 } \cdots  \sum_{\mathbf{j}_{s-1}\in\mathbb{N}_0 } 
p_{t}(\mathbf{n}_0)  f (\mathbf{j}_0,\mathbf{j}_1,\dots, \mathbf{j}_{s-1})
 \nonumber
 \\
 &\times&
\exp\left( {\mathbf{n}_0 \balpha^T} \right) \exp\left( \sum_{k=1}^{s} \mathbf{j}_{s-k} \left( \mathbf{H}^{(k-1)}(\bbeta)\right)^{T} \right)
 \exp\left( \mathbf{n}_{0} \left( \mathbf{H}^{(s)}(\bbeta)\right)^{T} \right),
 \nonumber
 \\
\label{eq:ape:inter4-multi-mgf}
\end{eqnarray}
which after rearranging the summands appropriately we get
\begin{eqnarray}
\mathcal{L}_{t,t+s} (\balpha,\bbeta)  &=&
\sum_{\mathbf{j}_0\in\mathbb{N}_0 } \sum_{\mathbf{j}_1\in\mathbb{N}_0 } \cdots  \sum_{\mathbf{j}_{s-1}\in\mathbb{N}_0 } 
 f (\mathbf{j}_0,\mathbf{j}_1,\dots, \mathbf{j}_{s-1})
\exp\left( \sum_{k=1}^{s} \mathbf{j}_{s-k} \left( \mathbf{H}^{(k-1)}(\bbeta)\right)^{T} \right)
 \nonumber
 \\
 &\times&
\sum_{\mathbf{n}_0\in\mathbb{N}_0 } p_{t}(\mathbf{n}_0)   \exp\left( \mathbf{n}_{0} \left( \balpha +\mathbf{H}^{(s)}(\bbeta)\right)^{T} \right)
\label{eq:ape:inter4-multi-mgf}
\end{eqnarray}
The summation over $\mathbf{j}_0,\mathbf{j}_1,\dots, \mathbf{j}_{s-1}$ turns out to be the moment generating function of the random vectors $\mathbf{J}^t,\mathbf{J}^{t+1},\dots,\mathbf{J}^{t+s}$, a function that will be denoted by $\mathcal{M}(\balpha_0,\balpha_1,\dots,\balpha_{s-1})$. Notice also that the summation over $\mathbf{n}_0$ results in the m.g.f.~of $\mathbf{N}^t$, which, as in preceding sections, has been denoted by $\mathcal{G}_t$. Using these conventions it is easy to see that, the m.g.f. of $(\mathbf{N}^t,\mathbf{N}^{t+s})$ is given by
\begin{eqnarray}
\mathcal{L}_{t,t+s} (\balpha,\bbeta)  &=&  \mathcal{G}_t \left( \balpha + \mathbf{H}^{(s)}(\bbeta) \right)
\mathcal{M}\left( \mathbf{H}^{(s-1)}(\bbeta),\mathbf{H}^{(s-2)}(\bbeta),\dots, \bbeta  \right).
\end{eqnarray}
Once we have computed the m.g.f.~it is possible to obtain the covariance matrix $\mathrm{Cov}(\mathbf{N}^t,\mathbf{N}^{t+s})$ by taking the second derivative of the above expression. We then obtain
\begin{eqnarray}
C_{i,j}(t,t+s) = \left( \mathrm{Cov}(\mathbf{N}^t,\mathbf{N}^{t+s})\right)_{i,j} = \frac{\partial^2  \mathcal{L}_{t,t+s} (\balpha,\bbeta)}{\partial \alpha_i \partial \beta_j}
 \Bigg|_{\balpha=\mathbf{0},\bbeta=\mathbf{0}}
 - \left(\bmu_t\right)_i\left(\bmu_{t+s}\right)_j
\end{eqnarray}
Performing some calculation we have
\begin{eqnarray}
C_{i,j}(t,t+s) &=& \frac{\partial}{\partial \alpha_i} \Bigg[ \sum_{k=1}^S \frac{\partial  \mathcal{G}_t \left( \mathbb{\bgamma} \right) }{\partial \gamma_k} \frac{\partial \left( \mathbf{H}^{(s)}(\bbeta)\right)_k}{\partial \beta_j}
\mathcal{M}\left( \mathbf{H}^{(s-1)}(\bbeta),\mathbf{H}^{(s-2)}(\bbeta),\dots, \bbeta  \right) 
\nonumber
\\
&+&
 \mathcal{G}_t \left( \balpha + \mathbf{H}^{(s)}(\bbeta) \right) 
 \sum_{m=1}^{s}\sum_{k=1}^S  \frac{\partial  \mathcal{M}\left( \bgamma_0,\bgamma_1,\dots,\bgamma_{s-1} \right) }{\partial (\bgamma_m)_k } 
 \frac{\partial {H}_k^{(m-1)}(\bbeta)}{\partial \beta_i}
 \bigg]
 \Bigg|_{\balpha=\mathbf{0},\bbeta=\mathbf{0}}
 \nonumber
 \\
 &-&  \left(\bmu_t\right)_i\left(\bmu_{t+s}\right)_j
 \nonumber
 \\
  &=& \sum_{k=1}^S \frac{\partial^2  \mathcal{G}_t \left( \mathbb{\bgamma} \right) }{\partial \gamma_i \partial \gamma_k} \frac{\partial \left( \mathbf{H}^{(s)}(\bbeta)\right)_k}{\partial \beta_j}
\mathcal{M}\left( \mathbf{H}^{(s-1)}(\bbeta),\mathbf{H}^{(s-2)}(\bbeta),\dots, \bbeta  \right)  \Bigg|_{\balpha=\mathbf{0},\bbeta=\mathbf{0}}
\nonumber
\\
&+&
\bigg( \frac{\partial \mathcal{G}_t \left( \balpha\right)}{\partial \alpha_i}  \bigg)
 \sum_{m=1}^{s}\sum_{k=1}^S  \frac{\partial  \mathcal{M}\left( \bgamma_0,\bgamma_1,\dots,\bgamma_{s-1} \right) }{\partial (\bgamma_m)_k } 
 \frac{\partial {H}_k^{(s-1-m)}(\bbeta)}{\partial \beta_i}
 \Bigg|_{\balpha=\mathbf{0},\bbeta=\mathbf{0}}
 \nonumber
 \\
 &-&  \left(\bmu_t\right)_i\left(\bmu_{t+s}\right)_j.
 \nonumber
\end{eqnarray}

Recalling that $  \mathcal{M} (\mathbf{0},\mathbf{0},\dots,\mathbf{0}) = 1$, $ \mathcal{G}_t(\mathbf{0})=1$, $\left( \partial  \mathcal{G}_t (\balpha)/\partial \alpha_j|_{\balpha = \mathbf{0}}\right) = (\bmu_t)_j$ and observing that
\begin{eqnarray}
\frac{\partial {H}_k^{(m-1)}(\bbeta)}{\partial \beta_j}\bigg|_{\bbeta=0} &=& \left( \mathbf{Q}^{m-1}\right)_{k,j}
\\
 \frac{\partial^2  \mathcal{G}_t \left( \mathbb{\bgamma} \right) }{\partial \gamma_i \partial \gamma_k} \Bigg|_{\bgamma=\mathbf{0}} &=& \left( \bSigma_t\right)_{i,k} +  \left(\bmu_t\right)_i\left(\bmu_{t}\right)_k
\\
\frac{\partial \mathcal{M}\left( \bgamma_0,\bgamma_1,\dots,\bgamma_{s-1} \right) }{\partial (\bgamma_m)_k } \bigg|_{\bbeta=0}   &=& (\bepsilon_{t+m})_k,
\end{eqnarray}
it is relatively easy to see that
\begin{eqnarray}
C_{i,j}(t,t+s)  &=& 
\sum_{k=1}^S   \left( \bSigma_t\right)_{i,k}   \left( \mathbf{Q}^{s}\right)_{k,j}   + \sum_{k=1}^S  \left(\bmu_t\right)_i\left(\bmu_{t}\right)_k   \left( \mathbf{Q}^{s}\right)_{k,j} 
\nonumber
\\
&+& (\bmu_t)_i \sum_{m=0}^{s-1}\sum_{k=1}^S  (\bepsilon_{t+m})_k \left( \mathbf{Q}^{s-1-m}\right)_{k,j} - \left(\bmu_t\right)_i\left(\bmu_{t+s}\right)_j
\nonumber
\\
&=&
\left( \bSigma_t  \mathbf{Q}^{s} \right)_{i,j} +  \left(\bmu_t\right)_i  \left( \bmu_t \mathbf{Q}^{s}\right)_{j}    +\sum_{m=1}^{s}  \left(\bmu_t \right)_i \left( \bepsilon_{t+m} \mathbf{Q}^{m-1}  \right)_{j} -  \left(\bmu_t\right)_i\left(\bmu_{t+s}\right)_j.
\nonumber
\\
\label{eq:ape:inter-cov}
\end{eqnarray}
Now, using the dynamics for $\bmu_t$ given by the equation
\[
\bmu_{t+1} = \bmu_t \mathbf{Q} + \bepsilon_t 
\]
it is easy to see that we can write $\bmu_{t+s}$ as
\begin{equation}
\bmu_{t+s} = \bmu_{t}\mathbf{Q}^s + \sum_{k=0}^{s-1} \bepsilon_{t+k} \mathbf{Q}^{s-1-k}.
\end{equation}
Using the above identity it is clear that equation~(\ref{eq:ape:inter-cov}) results in
\begin{equation}
\mathrm{Cov} (\mathbf{N}^t, \mathbf{N}^{t+s}) =  \bSigma_t  \mathbf{Q}^{s},
\end{equation}
which is the relationship anticipated in equation~(\ref{eq:covNtNt+s}).

\nocite{*}
\bibliography{OMC_ref}

\end{document}